\documentclass[seceqn,secthm]{elsart}
\usepackage{amssymb}\usepackage{amsmath}
\usepackage[all, knot]{xy}
\xyoption{arc} \xyoption{color}
\usepackage{mathrsfs}

\newcommand{\diag}{\operatorname{diag}}
\newtheorem{theorem}{Theorem}[section] % 1st argument is your name for it
\newtheorem{lemma}[theorem]{Lemma}     % 2nd argument is what is printed

\newtheorem{proposition}[theorem]{Proposition}
\newtheorem{remark}[theorem]{Remark}

\begin{document}

\begin{frontmatter}

\title{Complex scaling for the Dirichlet Laplacian in a domain with asymptotically cylindrical end\thanksref{AKA}}%
\thanks[AKA]{This work was funded by grant N108898 awarded by the Academy of Finland.}

\author{Victor Kalvin}
\ead{vkalvin\,@\,gmail.com} \journal{arXiv:0906.0601v2}
\date{}
\begin{abstract} We develop the complex scaling method for the Dirichlet Laplacian in a domain with asymptotically cylindrical end.  We define resonances as discrete eigenvalues of non-selfadjoint operators, obtained as deformations of the selfadjoint Dirichlet Laplacian $\Delta$ by means of the complex scaling. The resonances are identified with the poles of the  resolvent matrix elements $((\Delta-\mu)^{-1}F, G)$ meromorphic continuation in $\mu$    across the essential spectrum of $\Delta$, where $F$ and $G$ are elements of an explicitly given set of partial analytic vectors. It turns out that the Dirichlet Laplacian has no singular continuous spectrum, and its eigenvalues can accumulate only at threshold values of the spectral parameter.
\end{abstract}
\end{frontmatter}

%\begin{keywords}  asymptotically cylindrical ends, complex scaling, absolutely continuous spectrum, accumulations  of eigenvalues,  resonances
%\end{keywords}

%\tableofcontents

\section{Introduction} 
\begin{figure}[h]
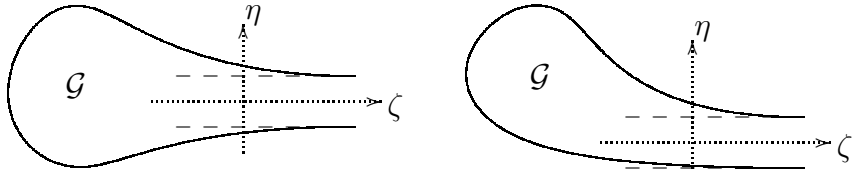

\[\xy (0,0)*{\xy0;/r.16pc/:
  (35,5); (35,-5) **\crv{(0,5)&(-20,30)&(-40,0)&(-20,-20)&(0,-5)};
 %(35,5); (35,-5) **\crv~pC{(0,5)&(-20,30)&(-40,0)&(-20,-30)&(0,-5)};
 {\ar@{.>} (-5,0)*{}; (40,0)*{}}; (43,-1)*{\zeta}; {\ar@{.>}
 (13,-10)*{}; (13,15)*{}};(15,17)*{\eta}; (-20,3)*{{\mathcal G}};
 (0,5)*{};(35,5)*{}**\dir{--};(0,-5)*{};(35,-5)*{}**\dir{--};\endxy};
%===================================================================
(60,0)*{\xy0;/r.16pc/: (35,5); (35,-5)
**\crv{(0,5)&(-10,35)&(-35,20)&(-30,0)&(0,-5)}; {\ar@{.>}
(-5,0)*{}; (40,0)*{}}; (43,-1)*{\zeta}; {\ar@{.>} (13,-5)*{};
(13,20)*{}}; (15,22)*{\eta};(-17,13)*{{\mathcal G}};
(0,5)*{};(35,5)*{}**\dir{--};(0,-5)*{};(35,-5)*{}**\dir{--};\endxy};
\endxy
\]
\caption{Domain ${\mathcal G}\subset \mathbb R^{n+1}$ with
asymptotically cylindrical end.}\label{fig1}
\end{figure}

We consider the  Dirichlet Laplacian in a domain $\mathcal G\subset\mathbb R^{n+1}$
with asymptotically cylindrical end. By the asymptotically
cylindrical end we mean any unbounded domain obtained as the
transformation of the semi-cylinder $(0,\infty)\times\Omega$ by some
suitable diffeomorphism $\varkappa$ that does not distort the space
at infinity; here $\Omega\subset\mathbb R^{n}$  is a bounded simply
connected domain. Outside of some compact set the domain $\mathcal
G$ coincides with its asymptotically cylindrical end, and the
boundary of $\mathcal G$ is smooth; see Fig.~\ref{fig1}.
 Equivalently, the asymptotically cylindrical end can be viewed as a
Riemannian manifold $((0,\infty)\times \Omega,\varkappa^*\mathsf e)$
isometrically embedded into the  space $\mathbb R^{n+1}$ with the
Euclidean metric $\mathsf e$. Then the pullback $\varkappa^*\mathsf
e$ of  $\mathsf e$ by the diffeomorphism $\varkappa$ converges in a
certain sense to $\mathsf e$  as the axial coordinate $x$ of the
semi-cylinder $(0,\infty)\times \Omega$ goes to infinity. Usually,
when studying the Laplacian in this kind of geometry, one imposes
some assumptions  on the rate of convergence of the metric
$\varkappa^*\mathsf e$ to its limit at infinity, e.g.~\cite{chr
ST,DES,DEM,Edward,FroHislop,Guillope,Melrose,MelroseScat,Mueller}.
Our goal is to study the case of arbitrarily slow convergence. With
this aim in mind we invoke the complex scaling method. This method
allows us to use an assumption on the analytic regularity of the
diffeomorphism $\varkappa$ with respect to the axial coordinate $x$
as a substitution for the assumptions on the rate of
 convergence of the metric $\varkappa^*\mathsf e$ at infinity.

The complex scaling method has a long tradition in the spectral and
scattering theories, and in numerical methods, see
e.g.~\cite{AC,BC,C,Cycon,DES,DEM,Hislop
Sigal,Hunziker,KalvinSiNum,MV1,MV2,Simon Reed iv,S,SZ,WZ} and
references therein. Nevertheless the complex scaling has not been
used in this geometric setting before. We first develop the complex
scaling method per se.  We introduce a deformation of the Dirichlet
Laplacian by means of the complex scaling, and  study the spectrum
of the deformed operator. Here we essentially rely on ideas of the
Aguilar-Balslev-Combes-Simon theory of
resonances~\cite{AC,BC,C,Cycon,Hislop Sigal,Simon Reed iv}, and most
of all on the approach to the complex scaling introduced
in~\cite{Hunziker}. As might be expected from known situations, the
deformation  by means of the complex scaling separates the essential
spectrum of the Dirichlet Laplacian $\Delta$ from the non-threshold
eigenvalues, rotating the rays of the essential spectrum about the
thresholds. For locating the essential spectrum we employ the theory
of elliptic boundary value
problems~\cite{KozlovMaz`ya,KozlovMazyaRossmann}. We characterize
the spectrum of the deformed Dirichlet Laplacian, establishing an
analog of the celebrated Aguilar-Balslev-Combes theorem. This
becomes possible due to the fact that on some sufficiently large
(dense) set of partial analytic vectors $\mathcal A$ the resolvent
matrix elements $((\Delta-\mu)^{-1}F, G)$ with $F,G\in\mathcal A$
are intimately connected with resolvent matrix elements of the
operator $\Delta$ deformed by means of the complex scaling. It turns
out that for all $F,G\in\mathcal A$ the function
$\mu\mapsto((\Delta-\mu)^{-1}F, G)$  can be meromorphically
continued from the physical sheet $\Bbb C\setminus (0,\infty)$
across the essential spectrum of $\Delta$ to a Riemann surface. This
implies that the Dirichlet Laplacian has no singular continuous
spectrum.   The non-real poles of the resolvent matrix elements are
identified with the resonances of the Dirichlet Laplacian, that are
introduced as the non-real eigenvalues of the deformed operator. The
real poles  correspond to the  eigenvalues of the Dirichlet
Laplacian. The eigenvalues embedded to the essential spectrum are
known to be very unstable, under rather weak perturbations they
become resonances~\cite{Aslan,Hislop Sigal}.

Under our assumptions on the domain $\mathcal G$, accumulations of
isolated and embedded eigenvalues of the Dirichlet Laplacian may
occur e.g.~\cite{Edward}. The Aguilar-Balslev-Combes theorem
implies that the non-threshold eigenvalues of the Dirichlet
Laplacian are of finite multiplicities, and  can accumulate only at
the thresholds. In a companion paper we will apply the complex
scaling in order to establish the following properties of the
Dirichlet Laplacian: a) the non-threshold eigenfunctions  are of
exponential decay at infinity; b) the eigenvalues are of finite
multiplicity; c) the eigenvalues can accumulate at the thresholds
only from below.
 Note that in~\cite{Edward} decay of eigenfunctions and eigenvalue accumulations of the Dirichlet and Neumann Laplacians are studied by another method and under an assumption on the form of the metric and on  the rate of its convergence at infinity.

We complete the introduction with a description of the structure of
this paper. In Section~\ref{s1} we formulate and discuss our
results. The subsequent sections are devoted to the proof.
Section~\ref{sec3} presents our approach to the complex scaling. In
Section~\ref{s5} we localize the essential spectrum of the Dirichlet
Laplacian deformed by means of the complex scaling. In
Section~\ref{s3} we
 consider a suitable set $\mathcal A$ of partial analytic vectors. Finally, in Section~\ref{s6} we construct the resolvent matrix elements meromorphic continuation and complete the proof of our results.

\section{Statement and discussion of results} \label{s1}
Let $(x,y)$ and $(\zeta,\eta)$ be two systems of the Cartesian
coordinates in $\mathbb R^{n+1}$, $n\geq 1$, such that
$x,\zeta\in\mathbb R$, while $y=(y_1,\dots,y_n)$ and
$\eta=(\eta_1,\dots,\eta_n)$ are in $\mathbb R^n$. Let
$\partial_x=\frac d {d x}$, $\partial_{y_m}=\frac d {d y_{m}}$,
  and $\partial_\zeta=\frac d {d \zeta}$, $\partial_{\eta_m}=\frac d {d \eta_{m}}$.

Consider a closed bounded simply connected domain
$\Omega\subset\mathbb R^n$, and the semi-cylinder $\Pi=\mathbb
R_+\times\Omega$, where $\mathbb R_+=\{x\in\mathbb R: x>0\}$. We say
that $\mathcal C\subset\mathbb R^{n+1}$ is an asymptotically
cylindrical end, if there exists a diffeomorphism
\begin{equation}\label{diff}
{\Pi}\ni(x,y)\mapsto \varkappa(x,y)=(\zeta,\eta)\in \mathcal C,
\end{equation}
such that the elements $\varkappa'_{\ell m}(x,\cdot)$ of its
Jacobian matrix $\varkappa'$ tend to the Kronecker delta
$\delta_{\ell m}$ in the space $C^\infty(\Omega)$  as $x\to+\infty$.

 Let
${\mathcal G}$ be a closed domain in $\mathbb R^{n+1}$ with smooth
boundary $\partial {\mathcal G}$. We suppose that the set
$\{(\zeta,\eta)\in {\mathcal G}: \zeta\leq 0\}$ is bounded, and the
set $\{(\zeta,\eta)\in {\mathcal G}:\zeta> 0\}$ is the
asymptotically cylindrical end $\mathcal C$, cf. Fig.~\ref{fig1}.
Introduce the notation
$\nabla_{\zeta\eta}=(\partial_\zeta,\partial_{\eta_1},\dots,\partial_{\eta_n})^\top$.
In the domain $\mathcal G$ we consider the Dirichlet Laplacian
$\Delta=-\nabla_{\zeta\eta}\cdot\nabla_{\zeta\eta}$ initially
defined on the core $C_0^\infty(\mathcal G)$. Here
$C_0^\infty(\mathcal G)$ is the set of all smooth
 compactly supported functions in $\mathcal
G$ satisfying the Dirichlet boundary condition
$u\upharpoonright_{\partial\mathcal G}=0$.

In this paper we use the complex scaling $x\mapsto x+\lambda v(x)$.
Here $\lambda$ is a scaling parameter, and $v(x)$ is  a smooth
function possessing the properties:
\begin{eqnarray}
& v(x)=0 \text{ for } x\leq R \text{ with a sufficiently large } R>0,\label{ab1}\\
& 0\leq v'(x)\leq 1 \text{ for all } x\in\mathbb R,\label{ab2}\\
& v'(x)=1 \text{ for } x\geq\widetilde{R}>R\label{ab3},
\end{eqnarray}
where $v'(x)=\partial_x v(x)$, and $\widetilde{R}$ is arbitrary. For
all real $\lambda\in(-1,1)$ the function $\mathbb R_+\ni x\mapsto
x+\lambda v(x)$ is invertible, and $\kappa_\lambda (x,y)=(x+\lambda
v(x),y)$ is a selfdiffeomorphism of the semi-cylinder $\Pi$.
Therefore
\begin{equation}\label{sc}
\vartheta_\lambda(\zeta,\eta)=\left\{
                      \begin{array}{lll}
                        \varkappa\circ\kappa_\lambda\circ\varkappa^{-1}(\zeta,\eta)&\text{ for }  & (\zeta,\eta)\in\mathcal C, \\
                       (\zeta,\eta) &\text{ for }& (\zeta,\eta)\in\mathcal G\setminus\mathcal C,
                      \end{array}
                    \right.
\end{equation}
is a selfdiffeomorphism of  the domain $\mathcal G$. In other words,
$\vartheta_\lambda$ with $\lambda\in(-1,1)$ is a scaling of the end
$\mathcal C$ along the curvilinear coordinate $x$.

 Let $(\vartheta_\lambda')^\top$ be the transpose of the Jacobian matrix $\vartheta_\lambda'$. Then
$\mathsf h_\lambda=(\vartheta_\lambda')^\top \vartheta_\lambda'$ is
the matrix coordinate representation of  a Riemannian metric
$\mathsf h_\lambda$ on $\mathcal G$, and
\begin{equation}\label{v2}
{^\lambda\!\Delta}=-\bigl({\det \mathsf h_\lambda
}\bigr)^{-1/2}\nabla_{\zeta\eta}\cdot\bigl({\det \mathsf h_\lambda
}\bigr)^{1/2}\mathsf h^{-1}_\lambda\nabla_{\zeta\eta}
\end{equation}
is the Laplace-Beltrami operator on the Riemannian manifold
$(\mathcal G,\mathsf h_\lambda)$.
  As the parameter $R$ in~\eqref{ab1} increases,  the equalities $\vartheta_\lambda(\zeta,\eta)=(\zeta,\eta)$,  $\mathsf
h_\lambda=\operatorname{Id}$, where $\operatorname{Id}$ is the
$(n+1)\times(n+1)$ identity matrix, and the equality
${^\lambda\!\Delta}=\Delta$,  become valid on a larger and larger
subset of $\mathcal G$.  In the case $\lambda=0$ the scaling is not
applied, and ${^0\!\Delta}\equiv\Delta$.

 In order to  consider  complex values of the scaling parameter $\lambda$, we  make  additional assumptions on the partial analytic regularity of the diffeomorphism $\varkappa$ in~\eqref{diff}:
\begin{enumerate}
\item[\it i.] the
function $\mathbb R_+\ni x\mapsto \varkappa(x,\cdot)\in
C^\infty(\Omega, \Bbb C^{n+1})$ has an analytic continuation from
$\mathbb R_+$ to some sector
    \begin{equation}\label{S}
    \mathbb S_\alpha=\{z\in\mathbb C: |\arg
    z|<\alpha<\pi/4\};
    \end{equation}
\item[\it ii.] the elements $\varkappa'_{\ell m}(z,\cdot)$ of the Jacobian matrix $\varkappa'$ tend to the
Kronecker delta $\delta_{\ell m}$ in the space $C^\infty(\Omega)$
uniformly in $z\in\mathbb S_\alpha$  as $z\to \infty$.
\end{enumerate}
Here $C^\infty(\Omega, \Bbb C^{n+1})$ is the space of smooth
functions acting from $\Omega$ to $\Bbb C^{n+1}$.

For instance, the assumptions~{\it i,ii}  are satisfied for the
following ends $\mathcal C\subset\mathbb R^2$:
$$
 \mathcal C=\{(\zeta,\eta)\in\mathbb R^2:(\zeta,\eta)=(x,y+\log(x+2)),\ x\in \mathbb R_+,\ y\in[0,1]\},
$$
$$
 \mathcal C=\Bigl\{(\zeta,\eta)\in \mathbb R^{2}: \zeta=\int_0^x\varphi(t)dt,\
 \eta=y\psi(x),\ x\in \mathbb R_+,\ y\in[0,1]\Bigr\},
$$
where as $\varphi(x)$ and $\psi(x)$ we can  take the functions $1$,
$1+e^{-x}$, $1+(x+1)^{-s}$ with $s>0$,  $1+1/\log(x+2)$,
$1+1/\log(1+\log(x+2))$, and so on.

It turns out that the assumptions {\it i, ii} together
with~\eqref{ab1} and \eqref{ab2}  lead to the analyticity of the
coefficients of the differential operator~\eqref{v2} with respect to
the scaling parameter $\lambda$ in the disc
\begin{equation}\label{disc}
\mathcal D_\alpha=\{\lambda\in\mathbb C: |\lambda|<\sin\alpha\}.
\end{equation}
We take the equality~\eqref{v2} as  definition of the differential
operator $^\lambda\!\Delta$ for all $\lambda\in\mathcal D_\alpha$.
The operator ${^\lambda\!\Delta}$ in $L^2(\mathcal G)$, initially
defined on the set $C_0^\infty(\mathcal G)$, is closable. Here
$L^2(\mathcal G)$ is the Hilbert space with the usual norm $ \|u;
L^2(\mathcal G)\|=\left(\int_{\mathcal G}
|u(\zeta,\eta)|^2\,d\zeta\,d\eta\right)^{1/2}$. The closure of
${^\lambda\!\Delta}$, denoted by the same symbol
${^\lambda\!\Delta}$, is  an unbounded operator  in $L^2(\mathcal
G)$, which is nonselfadjoint for $\lambda\neq 0$. However  the
Dirichlet Laplacian ${^0\!\Delta}\equiv\Delta$ is selfadjoint and
positive. We consider the operator ${^\lambda\!\Delta}$ with
$\lambda\in\mathcal D_\alpha$ as a deformation of $\Delta$  by means
of the complex scaling. The essential spectrum $\sigma_{ess}
({^\lambda\!\Delta})$ of  ${^\lambda\!\Delta}$ depends only on the
behaviour of the matrix $\mathsf h_\lambda$  outside  any compact
region of $\mathcal G$. In order to control $\sigma_{ess}
({^\lambda\!\Delta})$ we imposed the condition~\eqref{ab3}. Before
formulating our results, we introduce partial analytic vectors.

 Let $L^2(\Omega)$ be the space with the norm $(\int_\Omega |f(y)|^2\,dy)^{1/2}$. Consider the algebra  $\mathscr E$ of all entire functions $\mathbb C\ni z\mapsto f(z,\cdot)\in C^\infty(\Omega)$ with the following property:
in any sector $|\Im z|\leq (1-\epsilon) \Re z$ with $\epsilon>0$
the value $\|f(z,\cdot);L^2(\Omega)\|$ decays faster than any
inverse power of $\Re z$  as $\Re z\to+\infty$.  Examples of
functions $f\in \mathscr E$ are $f(z,y)=e^{-\gamma z^2}P(z,y)$,
where  $\gamma>0$ and $P(z,y)$ is an arbitrary polynomial  in $z$
with coefficients in $C^\infty(\Omega)$. We say that $F\in
L^2(\mathcal G)$ is a partial analytic vector, if $F\circ
\varkappa(x,y)=f(x,y)$  for some $f\in\mathscr E$ and all
$(x,y)\in\Pi$. The set of all partial analytic vectors is denoted by
$\mathcal A$. Later on we will show that $\mathcal A$ is  dense in
the space $L^2(\mathcal G)$.

\begin{theorem}\label{T1} Let $\nu_1<\nu_2<\dots$ be  the
 distinct eigenvalues of the selfadjoint Dirichlet Laplacian
 $\Delta_\Omega=-\partial_{\eta_1}^2-\cdots-\partial_{\eta_n}^2$ in the space $L^2(\Omega)$. Let the scaling parameter $\lambda$ be in the disc~\eqref{disc}, where $\alpha$ is the same as in the conditions {\it i, ii} on the diffeomorphism~$\varkappa$. Assume that the deformation $^\lambda\!\Delta$ of the Dirichlet Laplacian $\Delta$ is constructed with a smooth scaling function $v(x)$ satisfying the conditions~\eqref{ab1}--\eqref{ab3}. Then the following assertions are valid.

\begin{itemize}
\item[1.] The spectrum $\sigma(^\lambda\!\Delta)$ of the operator  $^\lambda\!\Delta$  is independent of the choice of  $v(x)$.

\item[2.]  $\mu$ is a point of the essential spectrum $\sigma_{ess}({^\lambda\!\Delta})$, if and only if
\begin{equation}\label{eq9}
\mu= \nu_j\text{ or }\arg(\mu-\nu_j)=-2\arg (1+\lambda)\text{ for
some }j\in\mathbb N.
\end{equation}

\item[3.]  $\sigma({^\lambda\!\Delta})=\sigma_{ess}({^\lambda\!\Delta})\cup\sigma_d({^\lambda\!\Delta})$, where $\sigma_d(^\lambda\!\Delta)$ is the discrete spectrum of ${^\lambda\!\Delta}$.

\item[4.] Let $\mu\in\sigma_d({^\lambda\!\Delta})$. As $\lambda$ changes continuously in the disk $\mathcal D_\alpha$, the point $\mu$ remains in $\sigma_d({^\lambda\!\Delta})$ as long as $\mu\in\mathbb C\setminus\sigma_{ess}({^\lambda\!\Delta})$.

\item[5.] Let $(\cdot,\cdot)$ stand for the inner product in $L^2(\mathcal G)$. For any $F,G\in \mathcal A$ the analytic  function
$\mathbb C\setminus\overline{\mathbb R_+}\ni\mu\mapsto\bigl(
(\Delta-\mu)^{-1}F,G\bigr)$ has a meromorphic continuation to the
set $\mathbb C\setminus\sigma_{ess}({^\lambda\!\Delta})$. Moreover,
$\mu$ is a pole of the meromorphic continuation with some  $F,G\in
\mathcal A$, if and only if $\mu\in \sigma_d({^\lambda\!\Delta})$.

\item[6.] A point  $\mu\in \mathbb R$, such that $\mu\neq\nu_j$ for all $j\in\mathbb N$, is an eigenvalue of the Dirichlet Laplacian $\Delta$, if and only if $\mu\in \sigma_d({^\lambda\!\Delta})$ with $\Im\lambda\neq 0$.

\item[7.] The Dirichlet Laplacian $\Delta$  has no singular
continuous spectrum.
\end{itemize}
\end{theorem}

 A similar result for the stationary
Schr\"{o}dinger operator in $\mathbb R^n$ is known as the
Aguilar-Balslev-Combes theorem, see e.g.~\cite{Hislop Sigal,Simon
Reed iv} and references therein.

\begin{figure}[h]
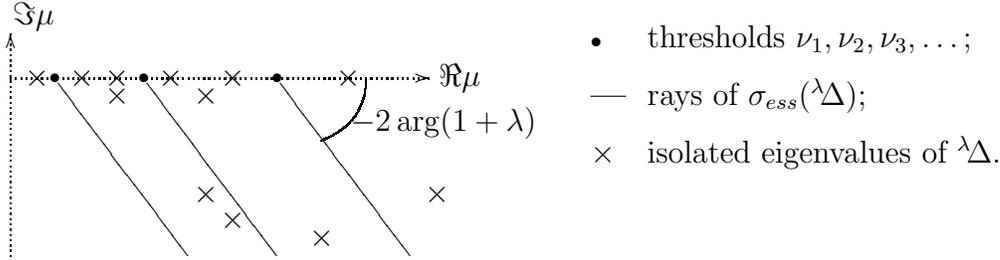

\[
\xy (0,0)*{\xy0;/r.28pc/:{\ar@{.>}(0,10);(50,10)*{\ \Re \mu}};
{\ar@{.>}(0,-10);(0,15)}; (2,17)*{\ \Im \mu};
(3,10)*{\times};(8,10)*{\times};(12,10)*{\times};(18,10)*{\times};(25,10)*{\times};(38,10)*{\times};{\ar@{-}(5,10)*{\scriptstyle\bullet};(20,-10)};{\ar@{-}(15,10)*{\scriptstyle\bullet};(30,-10)};
{(30,10)*{\scriptstyle\bullet};(45,-10)*{} **\dir{-}};
(40,10);(35,3)**\crv{(40,5)*{\quad\quad\quad\quad\quad
-2\arg(1+\lambda)}};
(35,-8)*{\times};(25,-6)*{\times};(22,-3)*{\times};(22,8)*{\times};(12,8)*{\times};(48,-3)*{\times}\endxy};
(70,0)*{\xy (20,20)*{
\begin{array}{ll}
     {\scriptstyle\bullet } & \text{ thresholds $\nu_1,\nu_2,\nu_3,\dots$;} \\
    \text{\bf---} & \text{ rays of $\sigma_{ess}({^\lambda\!\Delta})$;} \\
    {\times} &  \text{ isolated eigenvalues of $^\lambda\!\Delta$.}\\
\end{array}};
\endxy};
\endxy
\]
\caption{Spectral portrait of the deformed Dirichlet Laplacian
$^\lambda\!\Delta$, $\Im\lambda>0$.}\label{fig5}
\end{figure}

 The spectral portrait of the operator
${^\lambda \!\Delta}$ is depicted on Fig.~\ref{fig5}. We say  that
the numbers $\nu_1,\nu_2,\dots$ are threshold values of the spectral
parameter $\mu$, or thresholds for short.
 As the parameter
$\lambda$ varies, the ray $\arg(\mu-\nu_j)=-2\arg (1+\lambda)$ of
the essential spectrum
 $\sigma_{ess}({^\lambda\!\Delta})$ rotates about the threshold $\nu_j$,
and sweeps the sector $|\arg (\mu-\nu_j)|<2\alpha$. By the
assertion~{\it 4} the eigenvalues of ${^\lambda \!\Delta}$ outside
of the sector $|\arg (\mu-\nu_1)|<2\alpha$ do not change, hence they
are the discrete eigenvalues of the selfadjoint Dirichlet Laplacian
$\Delta$. All other discrete eigenvalues of ${^\lambda \!\Delta}$
belong to the sector $|\arg (\mu-\nu_1)|<2\alpha$. As $\lambda$
varies, they remain unchanged until they are covered by one of the
rotating rays of the essential spectrum. Conversely, new eigenvalues
can be uncovered by the rotating rays. In the case $\Im\lambda\geq
0$  (resp. $\Im\lambda\leq 0$) the operator ${^\lambda \!\Delta}$
 cannot have eigenvalues in the half-plane $\Im
\mu>0$ (resp. $\Im \mu<0$). Indeed, by the assertion~{\it 4} a
number $\mu$ with $\Im \mu>0$   is an eigenvalue of ${^\lambda
\!\Delta}$ with $\Im\lambda\geq 0$, if and only if $\mu$ is an
eigenvalue of $\Delta$, but  the Dirichlet Laplacian $\Delta$ cannot
have non-real eigenvalues as a selfadjoint operator. Further, by the
assertion~{\it 6} the real eigenvalues
$\mu\in\sigma_d({^\lambda\!\Delta})$ survive for $\lambda=0$: the
eigenvalues $\mu<\nu_1$ become the discrete eigenvalues of $\Delta$,
while the eigenvalues $\mu>\nu_1$ become the embedded non-threshold
eigenvalues of $\Delta$. In view of the fact that  any non-threshold
point $\mu$ can be separated from $\sigma_{ess}({^\lambda
\!\Delta})$ by a small variation of $\arg(1+\lambda)$,  the set
$\sigma_d({^\lambda \!\Delta})$ (and therefore the set of all
eigenvalues of $\Delta$) has no accumulation points, except possibly
for the thresholds $\nu_1,\nu_2,\dots$. Suitable examples of
accumulating eigenvalues can be found e.g. in~\cite{Edward}. By
definition, all discrete non-real eigenvalues of ${^\lambda
\!\Delta}$ are resonances of the Dirichlet Laplacian. By the
assertion~{\it 5} the resonances are characterized by the pair
$\{\Delta, \mathcal A\}$. They are identified with the complex poles
of the meromorphic continuation to a Riemann surface of all
resolvent matrix elements
$\mu\mapsto\bigl((\Delta-\mu)^{-1}F,G\bigr)$ with $F,G\in\mathcal
A$. Readers might have  noticed a certain analogy between the
situation we described above and the one known from the theory of
resonances for N-body problem e.g.~\cite{Hislop Sigal,Hunziker}.

\section{Complex scaling}\label{sec3}
In this section we show that the differential operator~\eqref{v2} is
well defined for complex values of the scaling parameter $\lambda$,
and obtain some estimates on its coefficients.

Consider the asymptotically cylindrical end $\mathcal C$ as a
Riemannian manifold endowed with the Euclidean metric $\mathsf e$.
We will use the coordinates $(\zeta,\eta)$ in $\mathcal G$ and
$(x,y)$ in $\Pi$, and identify the Riemannian metrics on $\Pi$ and
$\mathcal G$ with their matrix coordinate representations. Let
$\mathsf g=\varkappa^* \mathsf e$ be the pullback of $\mathsf e$ by
the diffeomorphism $\varkappa$ in~\eqref{diff}. Then the  matrix
$\mathsf g=[\mathsf g_{\ell m}]_{\ell ,m=1}^{n+1}$ is given by the
equality $\mathsf g=(\varkappa')^\top\varkappa'$,
 where $(\varkappa')^\top$ is the transpose of the Jacobian $\varkappa'$.
 Since the diffeomorphism $\varkappa$ satisfies the assumptions~{\it i,ii} of
Section~\ref{s1}, we conclude that  the metric matrix elements
 \begin{equation}\label{partan}
 \mathbb S_\alpha\ni z\mapsto  \mathsf
 g_{\ell m}(z,\cdot)\in C^\infty(\Omega)
 \end{equation}
 are analytic functions. Moreover, $\mathsf g_{\ell m}(z,\cdot)$ tends to the Kronecker delta $\delta_{\ell m}$ in the space $C^\infty(\Omega)$ uniformly in $z\in \mathbb S_\alpha$ as $z\to \infty$ or, equivalently, we have
 \begin{equation}\label{stab}
  \bigl|\partial_y^q(\mathsf g_{\ell m}(z,y)-
 \delta_{\ell m})\bigr |\leq C_q(|z|)\to 0  \text{ as } z\to \infty,\ z\in \mathbb S_\alpha,\ y\in\Omega,\ |q|\geq 0;
 \end{equation}
here
$\partial^q_{y}=\partial^{q_1}_{y_1}\partial^{q_2}_{y_2}\dots\partial^{q_n}_
{y_n}$ with a multiindex $q=(q_1,\dots,q_n)$, and  $|q|=\sum q_j$.

Consider the selfdiffeomorphism  $\kappa_\lambda(x,y)=(x+\lambda
v(x),y)$ of the semi-cylinder $\Pi$, where $\lambda\in(-1,1)$, and
$v(x)$ is a smooth scaling function with the
properties~\eqref{ab1}--\eqref{ab3}.  We define the metric $\mathsf
g_\lambda=\kappa_\lambda^*\mathsf g$ on $\Pi$ as the pullback of the
metric~$\mathsf g$ by $\kappa_\lambda$. As a result we get
 Riemannian manifolds $(\Pi, \mathsf g_\lambda)$ parametrized
by $\lambda\in(-1,1)$. For the matrix representation of  $\mathsf
g_\lambda$ we deduce the expression
 \begin{equation}\label{metric}
 \mathsf g_\lambda(x,y)=\diag\left\{1+\lambda
 v'(x),\operatorname{Id} \right\}\mathsf g(x+\lambda
 v(x),y)\diag\left\{1+\lambda v'(x),\operatorname{Id} \right\},
 \end{equation}
 where $\operatorname{Id}$
 stands for the $n\times n$-identity matrix, and $\diag\left\{1+\lambda v'(x),\operatorname{Id} \right\}$ is
 the Jacobian of $\kappa_\lambda$. The matrix $\mathsf g_\lambda$ is invertible for all $\lambda\in (-1,1)$. For brevity of notations we do not indicate the dependence of $v$, $\kappa_\lambda$, and $\mathsf g_\lambda$ on the large parameter $R$ in~\eqref{ab1}.

\begin{figure}
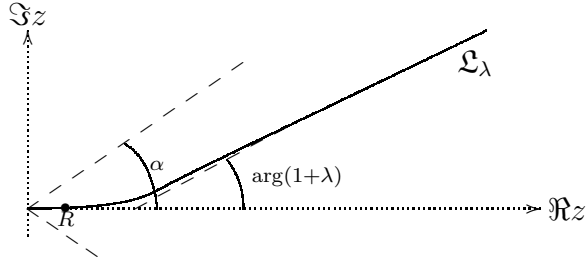
 \[ \xy0;/r.17pc/:
{\ar@{.>} (-12,0);(75,0)};(80,0)*{\Re z};
(-20,0);(20,27)**\dir{--};(-20,0);(-7,-9)**\dir{--};
(4,0);(-2,12)**\crv{(3,9)};(4,8)*{\scriptstyle\alpha};
{\ar@{.>} (-20,-5);(-20,33)}; (-20,36)*{ \Im z};%(40,40)*{\scriptstyle A_g (\lambda)};
(-20,0);(65,33)**\crv{(-2,0)&(5,4)&(15,9)&(36,19)};
(0,0);(28,15)**\dir{--};(-13,0)*{\scriptstyle\bullet};(-13,-2)*{\scriptstyle
R};
(17,9);(20,0)**\crv{(20,5)};(30,6)*{\scriptstyle\arg(1+\lambda)};
(63,27)*{ \mathfrak L_\lambda};
\endxy\]
\caption{The curve $\mathfrak L_\lambda$  for complex values of
$\lambda$.}\label{fig+}
\end{figure}

Now we wish to consider complex values of the scaling parameter
$\lambda$. We suppose that $\lambda$ is in the complex disc
$\mathcal D_\alpha$, where $\alpha$ is the same as in our
assumptions on the partial analytic regularity of the diffeomorphism
$\varkappa$; cf.~\eqref{S},~\eqref{disc}. The function $\mathbb
R_+\ni x\mapsto x+\lambda v(x)$ is invertible, and the curve
$\mathfrak L_\lambda=\{z\in\mathbb C: z=x+\lambda v(x),x>0 \}$ lies
in the sector $\mathbb S_\alpha$, see Fig.~\ref{fig+}.
We define the matrix $\mathsf g_\lambda$ for all non-real $\lambda$
in the disk by the equality~\eqref{metric}, where $\mathsf
g(x+\lambda v(x),y)$ stands for the value of the analytic in
$z\in\mathbb S_\alpha$ function $\mathsf g(z,y)$ at $z=x+\lambda
v(x)$.
 By
 analyticity in $\lambda$ we conclude that  ${\mathsf g}_\lambda$  is a
 complex symmetric matrix, the Schwartz reflection principle gives $\overline{{\mathsf g}_\lambda}= {{\mathsf g}_{\overline{\lambda}}}$, where the overline stands for the complex conjugation. If  $\lambda\in\mathcal D_\alpha$ is non-real, then the matrix  ${\mathsf g}_\lambda$ does not correspond to  a Riemannian metric.

Let us show that the matrix $\mathsf g_\lambda(x,y)$ is  invertible
for all $(x,y)\in\Pi$ and $\lambda\in\mathcal D_\alpha$. Evidently,
$v(x)=0$ and $\mathsf g^{-1}_\lambda(x,y)= \mathsf g^{-1}(x,y)$ for
all $x<R$. On the other hand, for all $x\geq R>0$ we have
$|x+\lambda v(x)|\geq  R$. Therefore for a sufficiently large $R>0$
the matrix $\mathsf g(x+\lambda v(x),y)$ with $x\geq R$ is little
different from the identity matrix due to~\eqref{stab}. This implies
invertibility of the matrix $\mathsf g(x+\lambda v(x),y)$ for all
$\lambda\in\mathcal D_\alpha$ and $(x,y)\in \Pi$. By analyticity the
equality
 \begin{equation}\label{metricinv}
 \mathsf g^{-1}_\lambda(x,y)=\diag\left\{\frac 1{1+\lambda
 v'(x)},\operatorname{Id} \right\}\mathsf
 g^{-1}(x+\lambda v(x),y)\diag\left\{\frac 1 {1+\lambda
 v'(x)},\operatorname{Id} \right\}
 \end{equation}
 extends  from real to all $\lambda\in\mathcal D_\alpha$.   Clearly,  the derivatives $\partial_x^p\partial_y^q \mathsf g^{-1}_\lambda$ are analytic functions of $\lambda\in\mathcal D_\alpha$.
On the next step we obtain some estimates on these derivatives.

Let $\partial_z=\frac 1 2(\partial_{\Re z}-i\partial_{\Im
 z})$ be the complex derivative.
Due to analyticity in $z$ we have $|\partial^p_z\partial_y^q\mathsf
g_{\ell m}(z,y)|\leq C_p\max _{\tilde z\in\mathcal
O(z)}|\partial_y^q\mathsf g_{\ell m}(\tilde z,y)|$, where $\mathcal
O(z)$ is a small circle centered at  $ z$. Therefore the
conditions~\eqref{partan} and~\eqref{stab} on the metric $\mathsf g$
lead to the uniform in  $y\in\Omega$ and $\{z\in\mathbb C:
|\arg(z-R)|<\alpha\}$  estimates
$$
\bigl|\partial_z^p\partial_y^q( \mathsf g_{\ell m}(z,y)-\delta_{\ell
m})\bigr|\leq \mathrm{c}_{pq}(R)\to 0 \text{ as } R\to+\infty,\quad
\ell,m=1,\dots,n+1.
$$
  This together with~\eqref{ab2} and~\eqref{metricinv} implies that
  \begin{equation}\label{lim}
\left \|\partial_x^p\partial_y^q\bigl(\mathsf g^{-1}_\lambda(x,y) -
\diag\bigl\{ (1+\lambda )^{-2},\operatorname{Id}
\bigr\}\bigr)\right\|\to 0 \text{ as } x\to+\infty,\  y\in\Omega,\
\lambda\in\mathcal D_\alpha,
\end{equation}
and  estimate
\begin{equation}\label{nstab}
\left \|\partial_x^p\partial_y^q\bigl(\mathsf g^{-1}_\lambda(x,y) -
\diag\bigl\{ (1+\lambda v'(x))^{-2},\operatorname{Id}
\bigr\}\bigr)\right\|\leq \mathrm{ C}_{pq}(R)
\end{equation}
holds uniformly in  $(x,y)\in[R,\infty)\times\Omega$ and
$\lambda\in\mathcal D_\alpha$, where $\|\cdot\|$ is the matrix norm
$\|A\|=\max_{\ell m}|a_{\ell m}|$, and $p+|q|\geq 0$. The constants
$\mathrm{C}_{pq}(R)$ in~\eqref{nstab} tend to zero as $R\to+\infty$;
recall that $v$ and $\mathsf g_\lambda$ both depend on the parameter
$R$.

 We define the complex scaling $\vartheta_\lambda$ for all $\lambda$ in the disk $\mathcal D_\alpha$ by the equality~\eqref{sc}, where $\varkappa\circ\kappa_\lambda(x,y)$  is  the value of the analytic in $z\in\mathbb S_\alpha$ function~$\varkappa(z,y)$ at the point $z=x+\lambda v(x)$. Consider the matrix $\mathsf h_\lambda=(\vartheta_\lambda')^\top \vartheta_\lambda'$. It is clear that for all  $(\zeta,\eta)\in\mathcal G
\setminus\mathcal C$ the matrix $\mathsf h_\lambda(\zeta,\eta)$
coincides with the $(n+1)\times(n+1)$-identity. For all real
$\lambda\in\mathcal D_\alpha$  the Riemannian geometry gives the
identity
\begin{equation}\label{rrr}
\mathsf g_\lambda(x,y)=\bigl(\varkappa'(x,y)\bigr)^\top\bigl(\mathsf
h_\lambda\circ\varkappa(x,y)\bigr)\varkappa'(x,y),\quad (x,y)\in\Pi,
\end{equation}
where  $\mathsf g_\lambda=\varkappa^*\mathsf h_\lambda$ is the
pullback of the corresponding metric $\mathsf h_\lambda$ on
$\mathcal C$ by the diffeomorphism $\varkappa$.  The Jacobian
$\varkappa'$ is an invertible and independent of $\lambda$ matrix,
and the matrix $\mathsf g_\lambda$ is invertible and analytic in
$\lambda\in\mathcal D_\alpha$. Therefore the matrix $\mathsf
h_\lambda$ is invertible and analytic in $\lambda\in\mathcal
D_\alpha$ due to~\eqref{rrr}.   By analyticity in $\lambda$ we
conclude that $\mathsf h_\lambda(\zeta,\eta)$ is a complex symmetric
matrix, the Schwartz reflection principle gives $\overline{\mathsf
h_\lambda}={\mathsf h_{\overline{\lambda}}}$.

Differentiating the equality~\eqref{rrr}, we see that  the
derivatives $\partial_\zeta^p\partial_\eta^q \mathsf h_\lambda $ and
$\partial_\zeta^p\partial_\eta^q \mathsf h^{-1}_\lambda$ are
analytic in $\lambda\in\mathcal D_\alpha$. Moreover,
from~\eqref{lim} and~\eqref{nstab} together with our assumptions  on
the diffeomorphism $\varkappa$ we have
$$
\left\|\partial^p_\zeta\partial^q_\eta\bigl(\mathsf
h^{-1}_\lambda(\zeta,\eta)- \diag\bigl\{ \bigl(1+\lambda
\bigr)^{-2},\operatorname{Id}\bigr\}\bigr)\right\|\to 0\text{ as
}\zeta\to+\infty,
$$
\begin{equation}\label{cont+}
{\begin{aligned}
&\left\|\partial^p_\zeta\partial^q_\eta\bigl(\mathsf
h^{-1}_\lambda(\zeta,\eta)- \diag\bigl\{ \bigl(1+\lambda
v'(\zeta,\eta)\bigr)^{-2},\operatorname{Id}
\bigr\}\bigr)\right\|\leq c(R),
\\
& \text{ where } p+|q|\leq 1,  \text{ and } \ c(R)\to 0 \text{ as }
R\to+\infty.
\end{aligned}}
\end{equation}
Here $(\zeta,\eta)\in\mathcal C$ and $v'(\zeta,\eta)\equiv
v'\circ\varkappa^{-1}(\zeta,\eta)$.   We extend  $v'$ from $\mathcal
C$ to $\mathcal G$ by zero, then the estimate~\eqref{cont+} extends
to all $(\zeta,\eta)\in\mathcal G$. The constant $c(R)$
in~\eqref{cont+} is independent of $\lambda\in\mathcal D_\alpha$ and
$(\zeta,\eta)\in\mathcal G$. Note that the matrix $\mathsf
h_\lambda$ with $\lambda\neq 0$ depends on the large parameter $R$,
however we do not indicate this in notations. Now we see that  the
differential operator~\eqref{v2} is well defined for all $\lambda$
in the disk $\mathcal D_\alpha$, and its coefficients are subjected
to the estimate~\eqref{cont+}.

\section{Localization of the essential spectrum}\label{s5}
 Introduce the
 the Sobolev space
$\lefteqn{\stackrel{\circ}{\phantom{\,\,>}}}H^2(\mathcal G)$ of
functions satisfying the Dirichlet boundary condition on
$\partial\mathcal G$ as the completion of the set
$C_0^\infty(\mathcal G)$ with respect to the norm
\begin{equation}\label{norm}
\|u;\lefteqn{\stackrel{\circ}{\phantom{\,\,>}}}H^2(\mathcal
G)\|=\|\Delta u; L^2(\mathcal G)\|+\|u;L^2(\mathcal G)\|.
\end{equation}

In this section we prove the following proposition.
\begin{proposition}\label{ess}
\begin{itemize}
\item[1.] The unbounded operator ${^\lambda\!\Delta}$ in $L^2(\mathcal G)$ with the domain $\lefteqn{\stackrel{\circ}{\phantom{\,\,>}}}H^2(\mathcal G)$ is closed.

\item[2.] The continuous operator ${^\lambda\!\Delta}-\mu:
\lefteqn{\stackrel{\circ}{\phantom{\,\,>}}}H^2(\mathcal G)\to
L^2(\mathcal G)$ is not  Fredholm, if and only if $\mu\in\mathbb C$
and $\lambda\in\mathcal D_\alpha$  meet the condition~\eqref{eq9}.
\end{itemize}
\end{proposition}
Recall that  $\mu$ is said to be a point of the essential spectrum
$\sigma_{ess}({^\lambda\!\Delta})$ of the closed unbounded operator
${^\lambda\!\Delta}$ in $L^2(\mathcal G)$ with the domain
$\lefteqn{\stackrel{\circ}{\phantom{\,\,>}}}H^2(\mathcal G)$, if the
continuous operator ${^\lambda\!\Delta-\mu}:
\lefteqn{\stackrel{\circ}{\phantom{\,\,>}}}H^2(\mathcal G)\to
L^2(\mathcal G)$ is not  Fredholm (a linear continuous operator
between two Banach spaces  is  Fredholm, if its kernel and cokernel
are finite dimensional, and its range is closed). Thus the
assertion~{\it 2} of Theorem~\ref{T1} is a direct consequence of
Proposition~\ref{ess}.

The proof of Proposition~\ref{ess} is essentially based on methods
of the theory  of elliptic non-homogeneous boundary value
problems~\cite{KozlovMaz`ya,KozlovMazyaRossmann,Lions Magenes}. The
proof  is preceded by the following lemma.
 \begin{lemma}\label{elliptic} The operator
 ${^\lambda\!\Delta}$ is strongly elliptic
for all $\lambda\in\mathcal D_\alpha$.
 \end{lemma}
 \begin{pf} Outside of the support of the scaling function $v$ the operator $^\lambda\!\Delta$ coincides with the strongly elliptic operator $\Delta$. Hence we only need to check the strong ellipticity of $^\lambda\!\Delta$ inside the image  of the set $[R,\infty)\times \Omega$ under the diffeomorphism  $\varkappa$.

 Consider the operator
 \begin{equation}\label{LB_l}
 {^\lambda\!\Delta_{\mathsf g}}=-\bigl(\det \mathsf
 g_\lambda \bigr)^{-1/2}\nabla_{xy}\cdot \bigl(\det \mathsf
 g_\lambda \bigr)^{1/2}\mathsf g^{-1}_\lambda
 \nabla_{xy},\quad \lambda\in\mathcal D_\alpha,
 \end{equation}
 where $\mathsf g_\lambda$ is the matrix~\eqref{metric} and $\nabla_{xy}=(\partial_x,\partial_{y_1}\dots
\partial_{y_n})^\top$. From~\eqref{rrr} and~\eqref{v2} it is easily  seen that $ {^\lambda\!\Delta_{\mathsf g}}$ is the operator $^\lambda\!\Delta$ written in the curvilinear coordinates $(x,y)$ inside the end $\mathcal C$; i.e. ${^\lambda\!\Delta} u=\bigl({^\lambda\!\Delta}_{\mathsf g}(u\circ\varkappa)\bigr)\circ\varkappa^{-1}$ for all $u\in C_0^\infty(\mathcal C)$. Since the strong ellipticity is preserved under the diffeomorphisms, it suffices to show that the operator $^\lambda\!\Delta_{\mathsf g}$  is strongly elliptic on $[R,\infty)\times \Omega\subset\Pi$.

By virtue of the bounds $|\lambda|<\sin\alpha$, $0<\alpha<\pi/4$,
and $0\leq v'(x)\leq 1$ we have
 \begin{equation}\label{symb}
 \Re\bigl(\xi\cdot \diag\bigl\{(1+\lambda v'(x))^{-2},\operatorname{Id}\bigr\}\xi\bigr)\geq (\cos 2\alpha) |\xi|^2/4
 \end{equation}
for all $x\in\mathbb R_+$ and $\xi\in\mathbb R^{n+1}$. Now we make
use of the estimates~\eqref{nstab} on  the matrix $\mathsf
g^{-1}_\lambda(x,y)$. Since $R$ is sufficiently large, the constant
$\mathrm{C}_{00}(R)$ in~\eqref{nstab} meets the estimate $(n+1)^{2}
\mathrm{C}_{00}(R)<(\cos 2\alpha)/4$. This together
with~\eqref{symb} implies the uniform in $\lambda\in\mathcal
D_\alpha$ and $(x,y)\in [R,\infty)\times \Omega$ estimates
 $$
 \begin{aligned}
 \Re \bigl(\xi\cdot & \mathsf g^{-1}_\lambda(x,y)  \xi\bigr)\geq \Re\bigl( \xi\cdot \diag\bigl\{(1+\lambda v'(x))^{-2},\operatorname{Id}\bigr\}\xi\bigr)\\
-&  \bigl|\xi\cdot\bigl(\mathsf g^{-1}_\lambda(x,y) - \diag\bigl\{
(1+\lambda )^{-2},\operatorname{Id} \bigr\}\bigr)\xi\bigr|
\\
&\geq \bigl( (\cos 2\alpha)/4
-(n+1)^2\mathrm{C}_{00}(R)\bigr)|\xi|^2
 \end{aligned}
 $$
  on the principal symbol of $\vphantom{D}^\lambda\!\Delta_{\mathsf g}$.
\end{pf}

\begin{pf*}{Proof of Proposition~\rm\ref{ess}}
We will rely on the following  lemma  due to Peetre, see
e.g.~\cite[Lemma~5.1]{Lions Magenes},
\cite[Lemma~3.4.1]{KozlovMazyaRossmann} or~\cite{Peetre}:
\begin{itemize}
\item[]{\it Let $\mathcal X,\mathcal Y$ and $\mathcal Z$ be Banach spaces, where $\mathcal X$ is compactly embedded into $\mathcal Z$. Furthemore, let $\mathcal L$ be a linear continuous operator from $\mathcal X$ to $\mathcal Y$. Then the next two assertions are equivalent: (i) the range of $\mathcal L$ is closed in $\mathcal Y$ and $\dim \ker \mathcal L<\infty$, (ii) there exists a constant $C$, such that
    \begin{equation}\label{1}
    \|u;{\mathcal X}\|\leq C(\|\mathcal L u;{\mathcal Y}\|+\|u;\mathcal Z\|)\quad \forall u\in \mathcal X.
    \end{equation}}
\end{itemize}
Below we assume that $\mu$ and $\lambda$ does not meet the
condition~\eqref{eq9} and establish the coercive estimate
 \begin{equation}\label{peetre}
\|u;H^2(\mathcal G)\|\leq C(\|({^\lambda\!\Delta}-\mu)u;
L^2(\mathcal G)\|+\|\mathsf w u; L^2(\mathcal G)\|)\quad \forall
u\in \lefteqn{\stackrel{\circ}{\phantom{\,\,>}}}H^2(\mathcal G)
\end{equation}
 for the operator ${^\lambda\!\Delta}-\mu: \lefteqn{\stackrel{\circ}{\phantom{\,\,>}}}H^2(\mathcal
G)\to L^2(\mathcal G)$. Here $\mathsf w\in C^\infty(\mathcal G) $ is
a positive rapidly decreasing at infinity  weight, such that the
embedding of
$\lefteqn{\stackrel{\circ}{\phantom{\,\,>}}}H^2(\mathcal G)$ into
the weighted space $L^2(\mathcal G;\mathsf w)$ with the norm
$\|\mathsf w \cdot; L^2(\mathcal G)\|$ is compact. The coercive
estimate~\eqref{peetre} is an estimate of type~\eqref{1}.

As is well-known, a strongly elliptic operator and the Dirichlet
boundary condition set up a regular elliptic boundary value problem,
e.g.~\cite{Lions Magenes}. Solutions of a regular elliptic boundary
value problem satisfy local coercive  estimates, e.g.~\cite{Lions
Magenes} or~\cite{KozlovMazyaRossmann}.  As a consequence of
Lemma~\ref{elliptic} we get the local coercive estimate
\begin{equation}\label{a p}
\|\rho_T u;H^2(\mathcal G)\|\leq
C(\|\varrho_T(^\lambda\!\Delta-\mu)u; L^2(\mathcal G)\|+\|\varrho_T
u; L^2(\mathcal G)\|)\quad\forall
u\in\lefteqn{\stackrel{\circ}{\phantom{\,\,>}}}H^2(\mathcal G).
\end{equation}
Here $\rho_T,\varrho_T\in C^\infty(\mathcal G)$ are compactly
supported cutoff functions, such that $\rho_T(\zeta,\eta)=1$ for
$|\zeta|<T+1$ and  $\varrho_T\rho_T=\rho_T$, where  $T$ is a  large
fixed number.

Let $\chi_T\in C^\infty(\mathcal G)$ be another cutoff function,
such that $\chi_T(\zeta,\eta)=1$ for $|\zeta|>T$ and
$\chi_T(\zeta,\eta)=0$ for $|\zeta|<T-1$. On the next step we
establish the estimate~\eqref{1} with $u$ replaced by  $\chi_T u$.
We will do it in the coordinates $(x,y)\in\Pi$.

Let $L^2(\Bbb R\times\Omega)$ be the space of functions in the
infinite cylinder $\Bbb R\times\Omega$ with the norm
$\bigl(\int_{\Bbb R}\|\mathsf u(x);
L^2(\Omega)\|^2\,dx\bigr)^{1/2}$. Introduce the Sobolev space
$\lefteqn{\stackrel{\circ}{\phantom{\,\,>}}}H^2(\Bbb R\times\Omega)$
of functions satisfying the Dirichlet boundary condition on $\mathbb
R\times\partial\Omega$ as the completion of the set $C_0^\infty(\Bbb
R\times\Omega)$ with respect to the norm
$$
\|\mathsf u;\lefteqn{\stackrel{\circ}{\phantom{\,\,>}}}H^2(\mathbb
R\times\Omega)\|=\Bigl(\sum_{ p+|q|\leq 2 }
 \|\partial_x^p\partial _y^q \mathsf u;
L^2(\Bbb R\times\Omega)\|^2\Bigr)^{1/2}.
$$
(Here $C_0^\infty(\Bbb R\times\Omega)$ is the set of all smooth
compactly supported functions $\mathsf u$ satisfying $\mathsf
u\upharpoonright_{\Bbb R\times\partial\Omega}=0$; recall that
$\Omega$ is closed.) Let $\mathsf u=(\chi_T u)\circ\varkappa$, where
$\varkappa$ is the diffeomorphism~\eqref{diff}. Due to our
assumptions on $\varkappa$ the estimates $0<\epsilon\leq\det
\varkappa'(x,y)\leq 1/\epsilon$ hold uniformly in $(x,y)\in\Pi$.
Hence for some independent of $u\in C_0^\infty(\mathcal G)$
constants $c_1$, $c_2$, and $c_3$ we have
\begin{equation}\label{www}
\begin{aligned}
&\|\chi_T u;\lefteqn{\stackrel{\circ}{\phantom{\,\,>}}}H^2(\mathcal
G)\|=\|\Delta (\chi_T u);L^2(\mathcal G)\|+\|\chi_T u; L^2(\mathcal
G)\|
\\
&\quad\leq c_1(\|{^0\!\Delta_{\mathsf g}}\mathsf u; L^2(\mathbb
R\times\Omega)\| +\|\mathsf u; L^2(\mathbb R\times\Omega)\|)\leq
c_2\|\mathsf u;\lefteqn{\stackrel{\circ}{\phantom{\,\,>}}} H^2(\Bbb
R\times\Omega)\|,
\\
&\|({^\lambda\!\Delta_{\mathsf g}-\mu})\mathsf u;  L^2(\mathbb
R\times\Omega)\|\leq c_3\|({^\lambda\!\Delta-\mu})\chi_T u;
L^2(\mathcal G)\|.
\end{aligned}
\end{equation}
Here the functions $\mathsf u$ and ${^\lambda\!\Delta_{\mathsf
g}}\mathsf u=({^\lambda\!\Delta} (\chi_T u))\circ\varkappa$ are
extended from $\Pi$ to the infinite cylinder $\mathbb R\times\Omega$
by zero, and $\|{^0\!\Delta_{\mathsf g}}\mathsf u; L^2(\mathbb
R\times\Omega)\|\leq C\|\mathsf
u;\lefteqn{\stackrel{\circ}{\phantom{\,\,>}}} H^2(\Bbb
R\times\Omega)\|$ because the coefficients of the Laplacian
${^0\!\Delta_{\mathsf g}}$ are bounded, cf.~\eqref{LB_l}
and~\eqref{nstab}. As $T$ is large, the function $\mathsf u$ is
supported in a small neighborhood of infinity. Due to the
stabilization condition~\eqref{lim} on $\mathsf g^{-1}_\lambda$ the
coefficients of the differential operator
${^\lambda\!\Delta}_{\mathsf g}- \Delta_\Omega+(1  +
\lambda)^{-2}\partial_x^2$ are small on the support of $\mathsf u$.
As a result we get the estimate
\begin{equation}\label{eq11}
\bigl\|\bigl({^\lambda\!\Delta}_{\mathsf g}- \Delta_\Omega+(1  +
\lambda)^{-2}\partial_x^2\bigr)\mathsf u;L^2(\mathbb
R\times\Omega)\bigr\| \leq \epsilon\|\mathsf u;
\lefteqn{\stackrel{\circ}{\phantom{\,\,>}}}H^2(\mathbb
R\times\Omega)\|,
\end{equation}
where $\epsilon$ is small and independent of $u\in
C^\infty_0(\mathcal G)$; moreover, $\epsilon\to 0$ as $T\to+\infty$.

Consider the continuous operator
\begin{equation}\label{mo}
\Delta_\Omega-(1+\lambda)^{-2}\partial_x^2-\mu:
\lefteqn{\stackrel{\circ}{\phantom{\,\,>}}}H^2(\mathbb
R\times\Omega)\to L^2(\mathbb R\times\Omega).
\end{equation}
Applying the Fourier transform $\mathscr F_{x\mapsto \tau}$ we pass
from the operator~\eqref{mo} to the Dirichlet Laplacian
$\Delta_\Omega+(1+\lambda)^{-2}\tau^2-\mu$ in $L^2(\Omega)$. Since
$\mu$ and $\lambda$ does not meet the condition~\eqref{eq9}, the
spectral parameter $\mu-(1+\lambda)^{-2}\tau^2$ is outside of the
spectrum $\{\nu_j\}_{j=1}^\infty$  of $\Delta_\Omega$ for all
$\tau\in\Bbb R$. Then a known argument, see e.g.~\cite[Theorem
5.2.2]{KozlovMazyaRossmann},~\cite[Theorem~2.4.1]{KozlovMaz`ya},
implies that the operator~\eqref{mo} realizes an isomorphism. In
particular the estimate
  \begin{equation}\label{coer}
\begin{aligned}
 \|\mathsf u; \lefteqn{\stackrel{\circ}{\phantom{\,\,>}}}H^{2}(\mathbb R\times\Omega)\| \leq  \mathrm{c} \bigl\|\bigl(\Delta_\Omega-(1+\lambda)^{-2}\partial_x^2 -\mu\bigr) \mathsf u; L^2(\mathbb R\times\Omega)\bigr\|
\end{aligned}
\end{equation}
is valid with an independent of $\mathsf u\in
\lefteqn{\stackrel{\circ}{\phantom{\,\,>}}}H^{2}(\mathbb
R\times\Omega)$ constant $\mathrm{c}$; in order to make the paper
selfcontained, we establish the estimate~\eqref{coer} in
Lemma~\ref{add} below the proof. As a consequence of~\eqref{coer}
and ~\eqref{eq11} we have
$$
\begin{aligned}
(1&-\epsilon\mathrm{c}) \|\mathsf u;
\lefteqn{\stackrel{\circ}{\phantom{\,\,>}}}H^{2}(\mathbb
R\times\Omega)\|\leq
\mathrm{c}\bigl\|\bigl(\Delta_\Omega-(1+\lambda)^{-2}\partial_x^2
-\mu\bigr) \mathsf u; L^2(\mathbb R\times\Omega)\bigr\|
\\
&-\mathrm{c}\bigl\|\bigl({^\lambda\!\Delta}_{\mathsf g}-
\Delta_\Omega+(1  + \lambda)^{-2}\partial_x^2\bigr)\mathsf
u;L^2(\mathbb R\times\Omega)\bigr\|\leq
\mathrm{c}\|({^\lambda\!\Delta_{\mathsf g}}-\mu)\mathsf u; L^2(\Bbb
R\times\Omega)\|.
\end{aligned}
$$
 If $T$ is sufficiently large, then $\epsilon \mathrm{c}<1$.
This together with~\eqref{www} gives
\begin{equation}\label{2}
\|\chi_T u;\lefteqn{\stackrel{\circ}{\phantom{\,\,>}}}H^2(\mathcal
G) \|\leq C\|({^\lambda\!\Delta-\mu})\chi_T u; L^2(\mathcal G)\|,
\end{equation}
where the constant $C=\mathrm{c}(1-\epsilon \mathrm{c})^{-1}c_2 c_3$
is independent of $u\in C_0^\infty(\mathcal G)$. By continuity the
estimate~\eqref{2} extends to all $u\in
\lefteqn{\stackrel{\circ}{\phantom{\,\,>}}}H^2(\mathcal G)$.

Now we combine~\eqref{2} with~\eqref{a p}, and arrive at the
estimates
\begin{equation}\label{3}
\begin{aligned}
\|u;\lefteqn{\stackrel{\circ}{\phantom{\,\,>}}}H^2(\mathcal G)\|\leq
\|\chi_T u;\lefteqn{\stackrel{\circ}{\phantom{\,\,>}}}H^2(\mathcal
G)\|+\|\rho_Tu;\lefteqn{\stackrel{\circ}{\phantom{\,\,>}}}H^2(\mathcal
G)\| \leq C(\|\chi_T({^\lambda\!\Delta-\mu}) u; L^2(\mathcal
G)\|\\+\|[{^\lambda\!\Delta-\mu},\chi_T] u; L^2(\mathcal G)\|
+\|\varrho({^\lambda\!\Delta-\mu}) u; L^2(\mathcal G)\|+\|\varrho u;
L^2(\mathcal G)\|)
\\
\leq C(\|({^\lambda\!\Delta-\mu}) u; L^2(\mathcal G)\|+\|\varrho u;
L^2(\mathcal G)\|).
\end{aligned}
\end{equation}
Here we used that $\rho_T=1$ on the support of the commutator
$[{^\lambda\!\Delta-\mu},\chi_T]$, and hence
$$
\|[{^\lambda\!\Delta-\mu},\chi_T] u; L^2(\mathcal G)\|\leq
C\|\rho_Tu;\lefteqn{\stackrel{\circ}{\phantom{\,\,>}}}H^2(\mathcal
G)\|.
$$

For an arbitrary positive weight $\mathsf w\in C^\infty(\mathcal G)$
we have $$ \|\varrho u; L^2(\mathcal G)\|\leq C\|\mathsf w u;
L^2(\mathcal G)\|$$ with an independent of
$u\in\lefteqn{\stackrel{\circ}{\phantom{\,\,>}}}H^2(\mathcal G)$
constant $C$. Thus the estimate~\eqref{1} is a direct consequence
of~\eqref{3}. By the Peetre's lemma we conclude that the range of
the continuous operator
${^\lambda\!\Delta}-\mu:\lefteqn{\stackrel{\circ}{\phantom{\,\,>}}}H^2(\mathcal
G)\to L^2(\mathcal G)$ is closed and the kernel is
finite-dimensional.
 The estimate~\eqref{1} with $\mathsf w\equiv 1$ shows that $\lefteqn{\stackrel{\circ}{\phantom{\,\,>}}}H^2(\mathcal
G)$ is the domain of the closed operator ${^\lambda\!\Delta}$ in
$L^2(\mathcal G)$, which proves the assertion~{\it 1}.

In order to see that the cokernel $\operatorname{coker}
({^\lambda\!\Delta}-\mu)=\ker ({^\lambda\!\Delta^*}-\overline\mu)$
of the continuous operator
${^\lambda\!\Delta}-\mu:\lefteqn{\stackrel{\circ}{\phantom{\,\,>}}}H^2(\mathcal
G)\to L^2(\mathcal G)$ is finite-dimensional  (if $\mu$ and
$\lambda$ does not meet the condition~\eqref{eq9}), we obtain the
coercive estimate
\begin{equation}\label{4}
\|u;\lefteqn{\stackrel{\circ}{\phantom{\,\,>}}}H^2(\mathcal G)\|\leq
C(\|({^\lambda\!\Delta^*}-\overline\mu)u; L^2(\mathcal
G)\|+\|\mathsf w u; L^2(\mathcal G)\|)
\end{equation}
for the adjoint ${^\lambda\!\Delta^*}$ of the unbounded operator
${^\lambda\!\Delta}$ in $L^2(\mathcal G)$, and apply the Peetre's
lemma. The proof of the estimate~\eqref{4} is similar to the proof
of~\eqref{1}, we omit it.

We have proved that the operator
${^\lambda\!\Delta}-\mu:\lefteqn{\stackrel{\circ}{\phantom{\,\,>}}}H^2(\mathcal
G)\to L^2(\mathcal G)$ is Fredholm, if $\lambda$ and $\mu$ do not
satisfy the condition~\eqref{eq9}. Now we assume that $\lambda$ and
$\mu$ satisfy the condition~\eqref{eq9} for some $j$, and show that
the operator
${^\lambda\!\Delta}-\mu:\lefteqn{\stackrel{\circ}{\phantom{\,\,>}}}H^2(\mathcal
G)\to L^2(\mathcal G)$ is not Fredholm.

Let $\chi$ be a smooth cutoff function on the real line, such that
$\chi(x)=1$ for $|x-3|\leq 1$ and $\chi(x)=0$ for $|x-3|\geq 2$.
Consider the functions
\begin{equation}\label{test}
\mathsf u_\ell(x,\mathrm
y)=\chi(x/\ell)\exp\bigl({i(1+\lambda)\sqrt{\mu-\nu_j}x}\bigr)\Phi(\mathrm
y),\quad (x,\mathrm y)\in \mathbb R\times\Omega,
\end{equation}
where $\Phi$ is an eigenfunction of $\Delta_\Omega$, corresponding
to the eigenvalue $\nu_j$. The exponent in~\eqref{test} is an
oscillating function of $x$.  Straightforward calculation shows that
\begin{equation}\label{result}
\bigl\|\bigl(\Delta_\Omega -(1+  \lambda)^{-2}\partial_x^2-\mu\bigr)
\mathsf u_\ell; L^2(\mathbb R\times\Omega)\bigr\|\leq
C,\quad\|\mathsf u_\ell;
\lefteqn{\stackrel{\circ}{\phantom{\,\,>}}}H^2(\Bbb
R\times\Omega)\|\to\infty
\end{equation}
 as $\ell\to +\infty$. Similarly to~\eqref{eq11} we conclude that
 \begin{equation}\label{5}
\bigl\|\bigl({^\lambda\!\Delta}_{\mathsf g}- \Delta_\Omega+(1  +
\lambda)^{-2}\partial_x^2\bigr)\mathsf u_\ell;L^2(\mathbb
R\times\Omega)\bigr\| \leq \epsilon_\ell\|\mathsf u_\ell;
\lefteqn{\stackrel{\circ}{\phantom{\,\,>}}}H^2(\Bbb
R\times\Omega)\|,
\end{equation}
where $\epsilon_\ell\to 0$ as $\ell\to+\infty$. Let the functions
$u_\ell=\mathsf u_\ell\circ\varkappa^{-1}$ be extended from
$\mathcal C$ to $\mathcal G$ by zero. If the operator
${^\lambda\!\Delta}-\mu:\lefteqn{\stackrel{\circ}{\phantom{\,\,>}}}H^2(\mathcal
G)\to L^2(\mathcal G)$ is Fredholm, then by the Peetre's lemma the
estimate~\eqref{peetre} holds with any weight $\mathsf w$, such that
$\lefteqn{\stackrel{\circ}{\phantom{\,\,>}}}H^2(\mathcal
G)\hookrightarrow L^2(\mathcal G;\mathsf w)$ is a compact embedding.
Without loss of generality we can assume that $\|\mathsf w u_\ell;
L^2(\mathcal G)\|\leq C$ for all $\ell\geq 1$. After the change of
variables $(\zeta,\eta)\mapsto(x,y)$ the estimate~\eqref{peetre}
implies
$$
\|\mathsf u_\ell;
\lefteqn{\stackrel{\circ}{\phantom{\,\,>}}}H^2(\Bbb
R\times\Omega)\|\leq C(\|({^\lambda\!\Delta_{\mathsf g}-\mu})\mathsf
u_\ell; L^2(\Bbb R\times\Omega)\|+1),
$$
where ${^\lambda\!\Delta_{\mathsf g}}\mathsf
u_\ell=({^\lambda\!\Delta}u_\ell)\circ\varkappa$ is extended from
$\Pi$ to $\Bbb R\times\Omega$ by zero. This together with~\eqref{5}
justifies the estimate
$$
\|\mathsf u_\ell;
\lefteqn{\stackrel{\circ}{\phantom{\,\,>}}}H^2(\Bbb
R\times\Omega)\|\leq C\bigl(\bigl\|\bigl(\Delta_\Omega -(1+
\lambda)^{-2}\partial_x^2-\mu\bigr) \mathsf u_\ell; L^2(\mathbb
R\times\Omega)\bigr\|+1\bigr),
$$
which contradicts~\eqref{result}.$\Box$
\end{pf*}

\begin{lemma}\label{add} Assume that  $|\lambda|<2^{-1/2}$ and $\mu\in\mathbb C$ do not meet the condition~\eqref{eq9}.
Then the estimate~\eqref{coer} holds with a constant $\mathrm{c}$
independent of $\mathsf u\in
\lefteqn{\stackrel{\circ}{\phantom{\,\,>}}}H^2(\mathbb
R\times\Omega)$.
\end{lemma}

\begin{pf} The differential operator
$\Delta_\Omega-(1+\lambda)^{-2}\partial^2_x-\mu$ with
$|\lambda|<2^{-1/2}$ is strongly elliptic, cf.~\eqref{symb}.
Therefore the local coercive estimate
\begin{equation}\label{cest}
\begin{aligned}
\|\varrho \mathsf u;
\lefteqn{\stackrel{\circ}{\phantom{\,\,>}}}H^2(\mathbb
R\times\Omega)\|^2\leq
c\Bigl(\bigl\|\varsigma\bigl(\Delta_\Omega-(1+\lambda)^{-2}\partial^2_x
-\mu\bigr) \mathsf u; L^2(\mathbb R\times\Omega)\bigr\|^2
\\+\|\varsigma \mathsf u; L^2(\mathbb R\times\Omega)\|^2\bigr)
\end{aligned}
\end{equation}
is valid,  where $\varrho$ and $\varsigma$ are smooth functions of
the variable $x$ with compact supports, and such that
$\varrho\not\equiv 0$, $\varrho\varsigma=\varrho$.

Introduce the Sobolev space
$\lefteqn{\stackrel{\circ}{\phantom{\,\,>}}}H^\ell(\Omega)$,
$\ell=0,1,2$, as the completion of the set $C_0^\infty(\Omega)$ with
respect to the norm $\sum_{|q|\leq\ell}\|\partial_y^q\Psi ;
L^2(\Omega)\|$. Note that
$\lefteqn{\stackrel{\circ}{\phantom{\,\,>}}}H^0(\Omega)\equiv
L^2(\Omega)$,
$\lefteqn{\stackrel{\circ}{\phantom{\,\,>}}}H^2(\Omega)$ is the
domain of the selfadjoint Dirichlet Laplacian $\Delta_\Omega$ in
$L^2(\Omega)$, and
$\lefteqn{\stackrel{\circ}{\phantom{\,\,>}}}H^1(\Omega)$ is the
domain of the quadratic form of $\Delta_\Omega$.
  We substitute $\mathsf u(x, y)=e^{i\tau x}\Psi(
y)$ with $\tau\in\mathbb R$ and $\Psi\in
\lefteqn{\stackrel{\circ}{\phantom{\,\,>}}}H^2(\Omega)$
into~\eqref{cest}. After simple manipulations we arrive at the
estimate
\begin{equation}\label{est0}
\begin{aligned}
\sum_{p=0}^2|\tau|^{2p}
\|\Psi;\lefteqn{\stackrel{\circ}{\phantom{\,\,>}}}
H^{2-p}(\Omega)\|^2 \leq
C\Bigl(\bigl\|\bigl(\Delta_\Omega+(1+\lambda)^{-2}\tau^2 -\mu\bigr)
\Psi; L^2(\Omega)\bigr\|^2\\+\|\Psi; L^2(\Omega)\|^2\Bigr),
\end{aligned}
\end{equation}
where the constant $C$ depends on $\varrho$ and $\varsigma$, but not
on $\tau$ or $\Psi$.

We also have the  elliptic coercive  estimate
\begin{equation}\label{+est0}
\begin{aligned}
 \|\Psi; \lefteqn{\stackrel{\circ}{\phantom{\,\,>}}}H^{2}(\Omega)\|^2 \leq  c\Bigl(\bigl\|\bigl(\Delta_\Omega+(1+\lambda)^{-2}\tau^2 -\mu\bigr) \Psi; L^2(\Omega)\bigr\|^2
+\|\Psi; L^2(\Omega)\|^2\Bigr)
\end{aligned}
\end{equation}
for the Dirichlet Laplacian in $\Omega$. Since $\lambda$ and $\mu$
do not meet the condition~\eqref{eq9}, the distance $d$ between the
ray $\{\mu-(1+\lambda)^{-2}\tau^2:\tau\in\mathbb R\}$ and the
spectrum $\{\nu_j\}_{j=1}^\infty$  of the selfadjoint operator
$\Delta_\Omega$ is positive. Therefore
$$
\|\Psi; L^2(\Omega)\|^2\leq
d^{-2}\bigl\|\bigl(\Delta_\Omega+(1+\lambda)^{-2}\tau^2 -\mu\bigr)
\Psi; L^2(\Omega)\bigr\|^2.
$$
The estimate~\eqref{+est0} takes the form
 \begin{equation}\label{est0-}
 \|\Psi; \lefteqn{\stackrel{\circ}{\phantom{\,\,>}}}H^{2}(\Omega)\|^2 \leq  (c+d^{-2})\bigl\|\bigl(\Delta_\Omega+(1+\lambda)^{-2}\tau^2 -\mu\bigr) \Psi; L^2(\Omega)\bigr\|^2.
\end{equation}

If $|\tau|> r$ with sufficiently large $r>0$, then the last term in
\eqref{est0} can be neglected. This together with~\eqref{est0-}
 justifies the estimate
\begin{equation}\label{estF}
\sum_{p=0}^2|\tau|^{2p}  \|\Psi;
\lefteqn{\stackrel{\circ}{\phantom{\,\,>}}}H^{2-p}(\Omega)\|^2 \leq
C \bigl\|\bigl(\Delta_\Omega+(1+\lambda)^{-2}\tau^2 -\mu\bigr) \Psi;
L^2(\Omega)\bigr\|^2
\end{equation}
 for all $\tau\in\mathbb R$ and some independent of $\Psi$ and $\tau$ constant $C$.

 Let $\Psi(\tau)$ be the Fourier transform $\int_\mathbb R e^{i\tau x} \mathsf u(x)\,dx$ of $\mathsf u\in \lefteqn{\stackrel{\circ}{\phantom{\,\,>}}} H^2(\mathbb R\times\Omega)$. Then $(-i\tau)^p \Psi(\tau)$ is the Fourier transform of $\partial_x^p \mathsf u(x)$. The Parseval equality gives
 $$
 \int_\mathbb R |\tau|^{2p} \|\Psi(\tau); \lefteqn{\stackrel{\circ}{\phantom{\,\,>}}}H^{2-p}(\Omega)\|^2\,d\tau=2\pi\int_\mathbb R \|\partial_x^p \mathsf u(x); \lefteqn{\stackrel{\circ}{\phantom{\,\,>}}}H^{2-p}(\Omega)\|^2\, dx.
 $$
 Integrating~\eqref{estF} with respect to $\tau\in\Bbb R$, we obtain the estimate~\eqref{coer}.
\end{pf}

\section{Partial analytic vectors}\label{s3}
Consider the set of partial analytic vectors $\mathcal A$,
introduced in Section~\ref{s1}.  Recall that  $F\in L^2(\mathcal G)$
is in the set $\mathcal A$, if $F\circ \varkappa(x,y)=f(x,y)$  for
some $f\in\mathscr E$ and all $(x,y)\in\Pi$. Here  $\mathscr E$ is
the algebra of all entire functions $\mathbb C\ni z\mapsto
f(z,\cdot)\in C^\infty(\Omega)$, such that in any sector $|\Im
z|\leq (1-\epsilon) \Re z$ with $\epsilon>0$  the value
$\|f(z,\cdot);L^2(\Omega)\|$ decays faster than any  inverse power
of $\Re z$  as $\Re z\to+\infty$.

For $F\in\mathcal A$ and $\lambda\in\mathbb C$  we define the
function  $F\circ\vartheta_\lambda$ in $\mathcal G$ by the
equalities $F\circ\vartheta_\lambda(\zeta,\eta)=F(\zeta,\eta)$ for
$(\zeta,\eta)\in\mathcal G\setminus\mathcal C$, and
\begin{equation}\label{eq}
F\circ\vartheta_\lambda\circ\varkappa(x,y)=f(x+\lambda v(x),y)\text{
for } (x,y)\in\Pi.
\end{equation}
Here $f(x+\lambda v(x),\cdot)$ is the value of the corresponding to
$F$ entire function $f\in\mathscr E$ at the point $z=x+\lambda
v(x)$.

\begin{proposition}\label{p1} Let $v$ be a smooth function satisfying~\eqref{ab1}--\eqref{ab3}. Then
\begin{itemize}
\item[1.] For any $F\in\mathcal A$, $\lambda\mapsto  F\circ\vartheta_\lambda$ is an $L^2(\mathcal G)$-valued analytic function in the disk $|\lambda|<2^{-1/2}$;
\item [2.] For any $\lambda$ in the disk $|\lambda|<2^{-1/2}$ the image
$\vartheta_\lambda[\mathcal A]=\{F\circ\vartheta_\lambda:
F\in\mathcal A\}$ of $\mathcal A$ under $\vartheta_\lambda$ is dense
in $L^2(\mathcal G)$.
\end{itemize}
\end{proposition}
\begin{pf}  In essence, this proposition is based on~\cite[Theorem 3]{Hunziker}.

Since $F\circ\vartheta_\lambda=F$ on $\mathcal G\setminus\mathcal
C$, it suffices to show that {\rm 1)} for any $F\in\mathcal A$,
$\lambda\mapsto  F\circ\vartheta_\lambda$ is an $L^2(\mathcal
C)$-valued analytic function in the disk $|\lambda|<2^{-1/2}$; {\rm
2)} for any $\lambda$ in this disk the image
$\vartheta_\lambda[\mathcal A]=\{F\circ\vartheta_\lambda:
F\in\mathcal A\}$ of $\mathcal A$ under $\vartheta_\lambda$ is dense
in $L^2(\mathcal C)$. Here the norm
\begin{equation}\label{n1}
\left(\int_{\mathcal C} | F(\zeta,\eta)|^2\,d\zeta\,
d\eta\right)^{1/2}=\left(\int_{\Pi}|
F\circ\varkappa(x,y)|^2\det\varkappa'(x,y)\,dx\,dy\right)^{1/2}
\end{equation}
in $L^2(\mathcal C)$ is equivalent to the norm
\begin{equation}\label{n2}
\left(\int_{\mathbb R_+} \| F\circ\varkappa(x,\cdot);
L^2(\Omega)\|^2\,dx\right)^{1/2}
\end{equation}
in $L^2(\Pi)$, because  $0<\epsilon<\det \varkappa'(x,y)<1/\epsilon$
uniformly in $(x,y)\in\Pi$ due to our assumptions on the
diffeomorphism $\varkappa$.

{\it 1.}   We set $z=x+\lambda v(x)$. Then $|\Re z|^2-|\Im z|^2\geq
x^2/2-|\lambda|^2|v(x)|^2$. Since $v(x)<x$, for all $\lambda$ in the
disk $|\lambda|\leq \sqrt{1/2-\epsilon}$, we get
$$
|\Re z|^2-|\Im z|^2\geq \epsilon x^2\geq \epsilon |v(x)|^2\geq
2\epsilon |\Im z|^2.
$$
Therefore $|\Im z|\leq (1+2\epsilon)^{-1/2}\Re z$. On the other hand
$\Re z\geq (1-2^{-1/2})x$. By the equality~\eqref{eq} with
$f\in\mathscr E$, combined with the definition of the
algebra~$\mathscr E$, we conclude that the value
$\|F\circ\vartheta_\lambda\circ\varkappa(x,\cdot); L^2(\Omega)\|$
decreases faster than any inverse power of $x$ as $x\to+\infty$,
uniformly in $\lambda$ with $|\lambda|\leq\sqrt{1/2-\epsilon}$. It
remains to note that $(
F\circ\vartheta_\lambda\circ\varkappa,G)_{\Pi}$ is analytic in
$|\lambda|<2^{-1/2}$ for any $ G \in L^2(\Pi)$, where
$(\cdot,\cdot)_{\Pi}$ stands for the inner product in $L^2(\Pi)$.
The assertion 1  is proven.

 {\it 2.} Given $h\in C^\infty_0(\Pi)$ we will construct
 a sequence $f_\ell\in\mathscr E$, such that  the
 function $\mathbb R_+\ni x\mapsto g_\ell(x)=f_\ell(x+\lambda v(x))\in C^\infty_0(\Omega)$ tends to $h$ in
 $L^2({\Pi})$ as $\ell\to +\infty$. Since the set $C_0^\infty(\Pi)$ is dense in $L^2(\Pi)$, and the norms~\eqref{n1} and~\eqref{n2} are equivalent, this will imply that the set $\vartheta_\lambda[\mathcal A]$ is dense in $L^2(\mathcal C)$.

 Namely, let
$$
f_\ell(z)=\sqrt{\frac {\ell} {\pi}}\int_{\mathbb R}  h(x)\exp[-\ell
(z-x-\lambda v(x))^2](1+\lambda v'(x))\,dx,\ \ell\geq 1,
$$
where $h(x)=0$ for $x\leq 0$. It is clear that $z\mapsto
 f_\ell(z)\in C_0^\infty(\Omega)$ is an entire function. Since $h$ is compactly supported, $z\mapsto \|f_\ell(z); L^2(\Omega)\|$  has the same falloff at infinity as $\exp(-\ell z^2)$, i.e. $f_\ell\in\mathscr E$.
In order to prove that $f_\ell$ tends to $h$ in $L^2(\Pi)$, we set
$$
s(x,\tilde x;\lambda)=x+\lambda v(x)-\tilde x-\lambda v(\tilde x).
$$
From the condition~\eqref{ab2} on the scaling function $v$ it
follows that for all $\lambda$ in the disk $|\lambda|\leq\sqrt{
1/2-\epsilon}$ we get $|\Re s|^2-|\Im s|^2\geq \epsilon|x-\tilde
x|^2$, and therefore
\begin{equation}\label{star}
|\exp(-s^2(x,\tilde x;\lambda))|\leq\exp(-\epsilon|x-\tilde x|^2).
\end{equation}
For all real $\lambda$ in the disk $|\lambda|\leq
\sqrt{1/2-\epsilon}$ we have
$$
\sqrt{\frac {\ell} {\pi}}\int_{\mathbb R} \exp[-\ell (x+\lambda
v(x)-\tilde x-\lambda v(\tilde x))^2](1+\lambda
v'(x))\,dx=\sqrt{\frac {\ell} {\pi}}\int_{\mathbb R}e^{-\ell
s^2}\,ds=1.
$$
Due to~\eqref{star}  these equalities extend by analyticity to the
disk $|\lambda|\leq \sqrt{1/2-\epsilon}$. Thus we established the
equality
\begin{equation}\label{eq14}
h(\tilde x)-g_\ell(\tilde x)=\sqrt{\frac {\ell} {\pi}}\int_{\mathbb
R}e^{-\ell s^2(x,\tilde x;\lambda)}(h(\tilde x)-h(x))(1+\lambda
v'(x))\, dx.
\end{equation}
This together with~\eqref{star} gives us the  estimate
\begin{equation}\label{eq13}
\|h(\tilde x)-g_\ell(\tilde x); L^2(\Omega)\|\leq C
\sqrt{\ell}\int_{\mathbb R}e^{-\ell\varepsilon (x-\tilde x)^2}
\|h(\tilde x)-h(x); L^2(\Omega)\|\,dx.
\end{equation}
It is known property of the Weierstra{\ss} singular
integral~\cite{Titchmarsh} that for all $\tilde x\in\mathbb R$
\begin{equation}\label{eq13+}
\sqrt{\ell}\int_{\mathbb R}e^{-\ell\varepsilon (x-\tilde
x)^2}\|h(\tilde x)-h(x); L^2(\Omega)\| \, dx \to 0\ \text{ as }\
\ell\to+\infty.
\end{equation}
Due to~\eqref{eq13} and~\eqref{eq13+}  $g_\ell$
 converges to $h$ in the norm of $L\sp 2({\Pi})$ as $\ell\to+\infty$.
\end{pf}

\section{Resolvent matrix elements meromorphic continuation}\label{s6}

In this section we construct the meromorphic continuation in $\mu$
of the  resolvent matrix elements $((\Delta-\mu)^{-1}F,G)$ with
$F,G\in\mathcal A$, and complete the proof of Theorem~\ref{T1}.

Let us first obtain a relation between the matrix elements
$((\Delta-\mu)^{-1}F,G)$ and some matrix elements of the resolvent
$({^\lambda}\Delta-\mu)^{-1}$ with a real $\lambda\in\mathcal
D_\alpha$.

For  $F,G\in L^2(\mathcal G)$ we define the sesquilinear form
$$
(F,G)_\lambda=\int_{\mathcal G}F\cdot \overline{G}\sqrt{\det\mathsf
h_\lambda}\,d\zeta\,d\eta,\quad \lambda\in\mathcal D_\alpha.
$$
As the parameter $R$ is sufficiently large, by the
estimate~\eqref{cont+} on $\mathsf h^{-1}_\lambda$ we have
$0<c_1\leq|\det\mathsf h_\lambda(\zeta,\eta)|\leq c_2$ uniformly in
$\lambda\in\mathcal D_\alpha$ and $(\zeta,\eta)\in\mathcal G$.
Therefore the form $(\cdot,\cdot)_\lambda$ in $L^2(\mathcal G)$ is
continuous and nondegenerate, i.e. $|(F,G)_\lambda|\leq C\|F;
L^2(\mathcal G)\|\|G; L^2(\mathcal G)\|$, and for any nonzero $F\in
L^2(\mathcal G)$  there exists $G\in L^2(\mathcal G)$, such that
$(F,G)_\lambda\neq 0$.

Assume that $\lambda\in\mathcal D_\alpha$ is real. Then the form
$(\cdot,\cdot)_\lambda$ is the inner product  induced on $\mathcal
G$ by the Riemannian metric $\mathsf h_\lambda$, and
$\sqrt{(\cdot,\cdot)_\lambda}$ is an equivalent norm in
$L^2(\mathcal G)$. The Laplace-Beltrami operator
${^\lambda\!\Delta}$  is a nonnegative selfadjoint operator in the
space $L^2(\mathcal G)$ endowed with  the norm
$\sqrt{(\cdot,\cdot)_\lambda}$, cf. Proposition~\ref{ess}.1. Thus
the resolvent $({^\lambda\!\Delta}-\mu)^{-1}$ with $\mu<0$ is
bounded. This allows to rewrite the equality $(\Delta
-\mu)u=\bigl(({^\lambda\!\Delta}-\mu)(u\circ\vartheta_\lambda)\bigr)\circ\vartheta_\lambda^{-1}$
with $u\in C_0^\infty(\mathcal G)$ in the form
\begin{equation}\label{++}
(\Delta
-\mu)^{-1}F=\bigl(({^\lambda\!\Delta}-\mu)^{-1}(F\circ\vartheta_\lambda)\bigr)\circ\vartheta_\lambda^{-1}.
\end{equation}
Here $F$ is in the set $\{F=(\Delta-\mu)u: u\in C^\infty_0(\mathcal
G)\}$. This set is dense in $L^2(\mathcal G)$, because
$C_0^\infty(\mathcal G)$ is dense in
$\lefteqn{\stackrel{\circ}{\phantom{\,\,>}}}H^2(\mathcal G)$, and
the operator
$\Delta-\mu:\lefteqn{\stackrel{\circ}{\phantom{\,\,>}}}H^2(\mathcal
G)\to L^2(\mathcal G)$ with $\mu<0$ realizes an isomorphism.

It is clear that
$(F\circ\vartheta_\lambda,F\circ\vartheta_\lambda)_\lambda=(F,F)$.
As a consequence, the (real) scaling $F\mapsto
F\circ\vartheta_\lambda$ realizes an isomorphism in $L^2(\mathcal
G)$, and the equality~\eqref{++} extends by continuity to all $F\in
L^2(\mathcal G)$. Taking the inner product  of the
equality~\eqref{++} with $G\in L^2(\mathcal G)$, and passing to the
variables $(\tilde\zeta,\tilde\eta)=\vartheta_\lambda(\zeta,\eta)$
in the right hand side, we obtain the relation
\begin{equation}\label{h6}
\bigl((\Delta-\mu)^{-1}F,G
\bigr)=\bigl(({^\lambda\!\Delta}-\mu)^{-1}(F\circ\vartheta_\lambda),G\circ\vartheta_{\overline{\lambda}}
\bigr)_\lambda, \quad\lambda\in \mathcal D_\alpha\cap\mathbb R,
\end{equation}
between the matrix elements of the resolvents $(\Delta-\mu)^{-1}$
and $({^\lambda\!\Delta}-\mu)^{-1}$.

We intend to implement the Aguilar-Balslev-Combes argument to the
equality~\eqref{h6}. In other words, for arbitrary $F,G\in\mathcal
A$ and for some fixed $\mu<0$ we will extend the equality~\eqref{h6}
by analyticity to all $\lambda$ in the disk $\mathcal D_\alpha$.
Then the right hand side of~\eqref{h6} will provide the left hand
side with a meromorphic continuation in $\mu$ across
$\sigma_{ess}({\Delta})$.

 By Proposition~\ref{p1} $F\circ\vartheta_\lambda$ and $G\circ\vartheta_\lambda$ are analytic $L^2(\mathcal G)$-valued functions of $\lambda\in\mathcal D_\alpha$.
 In order to implement the Aguilar-Balslev-Combes argument we need to show that the resolvent $({^\lambda\!\Delta}-\mu)^{-1}$ is an analytic function of  $\lambda$. With this aim in mind we prove the following proposition.
\begin{proposition}\label{LA}
The operator ${^\lambda\!\Delta}$ is
 m-sectorial with an independent of  $\lambda\in\mathcal D_\alpha$ sector.  Moreover,  $\mathcal
D_\alpha\ni\lambda\mapsto {^\lambda\!\Delta}$ is an analytic family
in the sense of Kato~\cite{Kato}.
\end{proposition}
\begin{remark} Here m-sectorial means that the numerical range and
the spectrum of ${^\lambda\!\Delta}$ are contained in some sector
$\{\mu\in\Bbb C: |\arg (\mu-a)|\leq b<\pi/2\}$.
\end{remark}
\begin{pf} The proof consists of
two steps. On the first step we introduce an operator
$[\vphantom{D}^\lambda\!\Delta]$ such that the difference
$\vphantom{D}^\lambda\!\Delta-[\vphantom{D}^\lambda\!\Delta]$ is a
first order differential operator. We estimate the numerical range
of $[{^\lambda\!\Delta}]$, and show that the operator
$[{^\lambda\!\Delta}]$ is m-sectorial. On the second step we prove
that ${^\lambda\!\Delta}$ is m-sectorial.

 {\it Step one.} We introduce the operator
\begin{equation}\label{op}
[{^\lambda\!\Delta}]=-\nabla_{\zeta\eta}\cdot \mathsf h^{-1}_\lambda
\nabla_{\zeta\eta}
\end{equation}
in the domain $\mathcal G$. From~\eqref{v2} and~\eqref{op} it is
clearly seen that
$\vphantom{D}^\lambda\!\Delta-[\vphantom{D}^\lambda\!\Delta]$ is a
first order operator. It is not hard to adapt the proof  of
Proposition~\ref{ess} for the operator $[{^\lambda\!\Delta}]$. It
turns out that
$\sigma_{ess}({^\lambda\!\Delta})=\sigma_{ess}([{^\lambda\!\Delta}])$,
and  the coercive estimate~\eqref{peetre} with   $^\lambda\!\Delta$
replaced by $[^\lambda\!\Delta]$ remains valid for all
$\mu\notin\sigma_{ess}({^\lambda\!\Delta})$. Therefore we can
consider $[{^\lambda\!\Delta}]$ as a closed unbounded operator in
$L^2(\mathcal G)$  with the domain
$\lefteqn{\stackrel{\circ}{\phantom{\,\,>}}}H^2(\mathcal G)$.

  In order to estimate the numerical range of  $[{^\lambda\!\Delta}]$, we
 establish
 the estimate
\begin{equation}
\label{eq7} \bigl|\arg \bigl([\vphantom{D}^\lambda\!\Delta] u,
u\bigr)\bigr|\leq 2\alpha+\sigma<\pi/2,\quad \forall u\in
\lefteqn{\stackrel{\circ}{\phantom{\,\,>}}}H^2(\mathcal G),
\end{equation}
for  all $\lambda\in\mathcal D_\alpha$; here $(\cdot,\cdot)$ is the
inner product in $L^2(\mathcal G)$. It is clear that
\begin{equation}\label{eq11+}
\begin{aligned}
&\bigl([{^\lambda\!\Delta}] u, u\bigr)=\int_{\mathcal
G}\bigl\langle\mathsf h_\lambda^{-1} \nabla_{\zeta \eta}u,
\nabla_{\zeta \eta}u\bigr\rangle\,d\zeta\,d\eta,
 \end{aligned}
\end{equation}
where $\langle\cdot,\cdot\rangle$ is the Hermitian inner product in
$\mathbb C^{n+1}$.

Let us estimate the numerical range of the matrix ${\mathsf
h}^{-1}_\lambda(\zeta,\eta)$. We shall rely on the
estimate~\eqref{cont+}. Let $Z\in\mathbb C^{n+1}$. Observe that by
virtue of the inequalities $0\leq v'(\zeta,\eta)\leq 1$ and
$|\lambda|<\sin\alpha<2^{-1/2}$ we have
 \begin{equation}\label{start}
 \begin{aligned}
& \Bigl|Z\cdot \diag\bigl\{\bigl(1+\lambda
v'(\zeta,\eta)\bigr)^{-2},\operatorname{Id}\bigr\}\overline
{Z}\Bigr|\geq |Z|^2/4,
 \\
& \Bigl|\arg \Bigl(Z\cdot \diag\bigl\{\bigl(1+\lambda
v'(\zeta,\eta)\bigr)^{-2},\operatorname{Id}\bigr\}\overline
{Z}\Bigr)\Bigr|<2\alpha.
\end{aligned}
 \end{equation}
Since $R$ is sufficiently large, the constant $c (R)$
in~\eqref{cont+} meets the estimate $4(n+1)^{2} c
(R)\leq\sin(\sigma/2)$ with some $\sigma\in(0,\pi/2-2\alpha)$.
Then~\eqref{cont+} together with~\eqref{start} gives
$$
\Bigl|\arg \Bigl(Z\cdot \mathsf h^{-1}_\lambda(\zeta,\eta)\overline
{Z}\Bigr)\Bigr|\leq 2\alpha+\sigma<\pi/2,\quad \lambda\in\mathcal
D_\alpha,\ (\zeta,\eta)\in\mathcal G.
$$
Taking into account~\eqref{eq11+} we arrive at~\eqref{eq7}.

To the operator $[{^\lambda\!\Delta}]$ with $\lambda\in\mathcal
D_\alpha$  we assign the sector $|\arg \mu|\leq 2\alpha+\sigma$ of
angle less than $\pi$. Due to~\eqref{eq7} the numerical range of
$[{^\lambda\!\Delta}]$ is inside of its sector. Moreover, as we
already found the essential spectrum
$\sigma_{ess}([{^\lambda\!\Delta}])=\sigma_{ess}({^\lambda\!\Delta})$,
we immediately see that it is  in the sector of
$[{^\lambda\!\Delta}]$, cf. Fig.~\ref{fig5}.

Let $\mu$ be outside of the sector of $[{^\lambda\!\Delta}]$. Then
the kernel of the operator $[{^\lambda\!\Delta}]-\mu$  is trivial.
In order to see that $[{^\lambda\!\Delta}]$ is m-sectorial, it
remains to show that the kernel of the adjoint operator
$([{^\lambda\!\Delta}]-\mu)^*$ is trivial.

Since the symmetric matrix $\mathsf h_\lambda$ meets the equality
$\overline{\mathsf h_\lambda}=\mathsf h_{\overline\lambda}$,  we
obtain the identity
$([{^\lambda\!\Delta}]u,v)=(u,[{^{\overline\lambda}\!\Delta}]v)$ for
all $u,v \in \lefteqn{\stackrel{\circ}{\phantom{\,\,>}}}H^2(\mathcal
G)$, cf.~\eqref{eq11+}. Hence the closed densely defined operator
$[{^{\overline\lambda}\!\Delta}]$ is adjoint to
$[{^\lambda\!\Delta}]$. Observe that $\overline{\mu}$ is outside of
the sector of $[{^{\overline\lambda}\!\Delta}]$, provided that $\mu$
is outside of the sector of $[{^\lambda\!\Delta}]$. Therefore the
kernel  of $[{^{\overline\lambda}\!\Delta}]-\overline{\mu}$ is
trivial for all $\mu$  outside of the sector of
$[{^\lambda\!\Delta}]$, and the operator $[{^\lambda\!\Delta}]$ is
m-sectorial.

{\it Step two.} Let us show that the  operator
$\vphantom{D}^\lambda\!\Delta-[\vphantom{D}^\lambda\!\Delta]$ has an
arbitrarily small and uniform in $\lambda\in\mathcal D_\alpha$
relative bound with respect to $[\vphantom{D}^\lambda\!\Delta]$ in
the operator sense.

 From the estimate~\eqref{cont+} on $\mathsf h^{-1}_\lambda$ and the definitions~\eqref{v2} and~\eqref{op} of the operators ${^\lambda\!\Delta}$ and $[{^\lambda\!\Delta}]$, we see that  coefficients $a^\lambda_{pq}(\zeta,\eta)$ of the first order
differential operator
${^\lambda\!\Delta}-[{^\lambda\!\Delta}]=\sum_{p+|q|=1}
a^\lambda_{pq}\partial_\zeta^p\partial_\eta^q$ are bounded uniformly
in $(\zeta,\eta)\in\mathcal G$ and $\lambda\in\mathcal D_\alpha$. By
a standard argument, based on the integration by parts,  we conclude
that the operator ${^\lambda\!\Delta}-[{^\lambda\!\Delta}]$ is
uniformly bounded relative to $\Delta$ with an arbitrarily small
relative bound; i.e. for all $\lambda\in\mathcal D_\alpha$ the
estimate
\begin{equation}\label{e}
\|({^\lambda\!\Delta}-[{^\lambda\!\Delta}])u; L^2(\mathcal G)\|\leq
\delta\|\Delta u; L^2(\mathcal G)\|+c(\delta)\|u; L^2(\mathcal
G)\|,\quad \forall u\in
\lefteqn{\stackrel{\circ}{\phantom{\,\,>}}}H^2(\mathcal G),
\end{equation}
holds with an arbitrarily small $\delta>0$ and some $c(\delta)$,
which depends only on $\delta$.

In Lemma~\ref{LA+} below we establish the coercive  estimate
\begin{equation}\label{a priori}
\|u;\lefteqn{\stackrel{\circ}{\phantom{\,\,>}}}H^2(\mathcal G)\|\leq
C(\|[{^\lambda\! \Delta}] u; L^2(\mathcal G)\|+\|u; L^2(\mathcal
G)\|),\quad \forall  u\in
\lefteqn{\stackrel{\circ}{\phantom{\,\,>}}}H^2(\mathcal G).
\end{equation}
 This estimate is far not
that sharp as the estimate~\eqref{peetre}, however it holds
uniformly in $\lambda\in\mathcal D_\alpha$. As a consequence
of~\eqref{norm},~\eqref{e}, and~\eqref{a priori} we obtain the
estimate
\begin{equation}\label{ee}
\|({^\lambda\!\Delta}-[{^\lambda\!\Delta}])u; L^2(\mathcal G)\|\leq
\delta\| [{^\lambda\!\Delta}]u; L^2(\mathcal G)\|+C(\delta)\|u;
L^2(\mathcal G)\|,
\end{equation}
where $\delta>0$ is arbitrarily small, and $C(\delta)$ depends only
on $\delta$.

The operator ${^\lambda\!\Delta}$ is m-sectorial as a perturbation
with an arbitrarily small relative bound of the m-sectorial operator
 $[{^\lambda\!\Delta}]$, e.g.~\cite{Kato,Simon Reed iv}. Moreover, the sector of ${^\lambda\!\Delta}$
is independent of  $\lambda$, because the estimate~\eqref{ee} holds
uniformly in $\lambda$, and the sector of $[{^\lambda\!\Delta}]$ is
independent of $\lambda$.

The coefficients of the m-sectorial operator ${^\lambda\!\Delta}$
are analytic with respect to $\lambda\in\mathcal D_\alpha$.
Therefore $\mathcal D_\alpha\ni\lambda\mapsto
\vphantom{D}^\lambda\!\Delta$ is an analytic family of type A; i.e.
for every $\lambda\in\mathcal D_\alpha$ the resolvent set of the
closed unbounded operator ${^\lambda\!\Delta}$ in $L^2(\mathcal G)$
is not empty, the domain
$\lefteqn{\stackrel{\circ}{\phantom{\,\,>}}}H^2(\mathcal G)$ of
${^\lambda\!\Delta}$ does not depend on $\lambda$, and for any
$u\in\lefteqn{\stackrel{\circ}{\phantom{\,\,>}}}H^2(\mathcal G)$ the
function $\mathcal D_\alpha\ni\lambda\mapsto
\vphantom{D}^\lambda\!\Delta u\in L^2(\mathcal G)$ is analytic. As
is known~\cite{Kato,Simon Reed iv}, an analytic family of type A is
also  analytic  in the sense of Kato.
\end{pf}

\begin{lemma}\label{LA+} The   coercive
estimate~\eqref{a priori} holds uniformly in $\lambda\in\mathcal
D_\alpha$.
\end{lemma}
\begin{pf} Consider the auxiliary operator
$$
\mathsf A_\lambda=-\partial_\zeta \bigl(1+\lambda
v'(\zeta,\eta)\bigr)^{-2}\partial_\zeta
-\partial_{\eta_1}^2-\cdots-\partial_{\eta_n}^2.
$$
 By the estimate~\eqref{cont+}  the
coefficients of the differential operator $\mathsf
A_\lambda-[^\lambda\!\Delta]$  can be made arbitrarily small
(uniformly in $\lambda\in\mathcal D_\alpha$) by taking a
sufficiently large $R$.
 Therefore for any $\epsilon>0$ there exists $R>0$, such that
\begin{equation}\label{h2}
\|\mathsf A_\lambda u-[{^\lambda\!\Delta}]u; L^2(\mathcal G)\|^2\leq
\epsilon\| u;
\lefteqn{\stackrel{\circ}{\phantom{\,\,>}}}H^2(\mathcal G )\|^2.
\end{equation}
 Below we will prove the uniform in $\lambda$ estimate
\begin{equation}\label{h1}
\|\Delta u; L^2(\mathcal G)\|^2\leq C(\|\mathsf A_\lambda u;
L^2(\mathcal G)\|^2+\|u; L^2(\mathcal G)\|^2)
\end{equation}
with some independent of $R$ and $u$ constant $C$. Then the coercive
stimate~\eqref{a priori} will follow from~\eqref{h2},~\eqref{h1},
and~\eqref{norm}.

Let us prove~\eqref{h1} for $u\in C_0^\infty(\mathcal G)$. From the
inequalities $0\leq v'(\zeta,\eta)\leq 1$, $|\lambda|<\sin\alpha$,
and $\alpha<\pi/4$, it is easily seen that
$$
\Re (\mathsf A_\lambda u,u)=\Re\bigl((1+\lambda
v')^{-2}\partial_\zeta u,\partial_\zeta
u\bigr)+\sum\bigl(\partial_{\eta_j}u,\partial_{\eta_j}u\bigr)\geq
(\Delta u,u)/c
$$
with $1/c=(\cos 2\alpha)/4>0$. As a consequence we have
\begin{equation}\label{h3}
\begin{aligned}
\|\Delta u; L^2(\mathcal
G)\|^2=(\Delta\nabla_{\zeta\eta}u,\nabla_{\zeta\eta}u)\leq
c|(\mathsf A_\lambda\nabla_{\zeta\eta} u,\nabla_{\zeta\eta}u)|
\\
 \leq c |(\nabla_{\zeta\eta}\mathsf
A_\lambda u,\nabla_{\zeta\eta}u)| +c|([\mathsf
A_\lambda,\nabla_{\zeta\eta}]u,\nabla_{\zeta\eta}u)|.
\end{aligned}
\end{equation}
 In the first term of the right hand
side of~\eqref{h3} we use that
\begin{equation}\label{n++}
 |(\nabla_{\zeta\eta}\mathsf
A_\lambda u,\nabla_{\zeta\eta}u)|= |(\mathsf A_\lambda u,\Delta
u)|\leq \delta \|\Delta u; L^2(\mathcal G)\|^2+\delta^{-1}\|\mathsf
A_\lambda u; L^2(\mathcal G)\|^2
\end{equation}
for an arbitrarily small $\delta>0$. The commutator $[\mathsf
A_\lambda,\nabla_{\zeta\eta}]$ in the second term is a second order
operator satisfying the estimate
$$
\|[\mathsf A_\lambda,\nabla_{\zeta\eta}]u; L^2(\mathcal G)\|\leq
C(\|\Delta u; L^2(\mathcal G)\|+\|\nabla_{\zeta\eta}u; L^2(\mathcal
G)\|).
$$
As a consequence, for an arbitrarily small $\delta>0$ we get
\begin{equation}\label{n+}
\begin{aligned}
&|([\mathsf A_\lambda,\nabla_{\zeta\eta}]
u,\nabla_{\zeta\eta}u)|\leq C(\|\Delta u; L^2(\mathcal
G)\|+\|\nabla_{\zeta\eta}u; L^2(\mathcal G)\|)
 \|\nabla_{\zeta\eta}u; L^2(\mathcal G)\|
 \\
 &\leq
\delta\|\Delta u; L^2(\mathcal G)\|^2+
c(\delta)(\nabla_{\zeta\eta}u,\nabla_{\zeta\eta}u) \leq 2\delta
\|\Delta u; L^2(\mathcal G)\|^2+ C(\delta)\|u;L^2(\mathcal G)\|^2.
\end{aligned}
\end{equation}
Now for all $u\in C_0^\infty(\mathcal G)$ the estimate~\eqref{h1}
is readily apparent from~\eqref{h3},~\eqref{n++}, and~\eqref{n+}. By
continuity it extends to all $u\in
\lefteqn{\stackrel{\circ}{\phantom{\,\,>}}}H^2(\mathcal G)$
\end{pf}

Now we are in position to prove Theorem~\ref{T1}.
\begin{pf*}{Proof of Theorem~\rm\ref{T1}} The assertion~{\it 2} is a direct consequence of Proposition~\ref{ess}.

{\it 3.} Let $\lambda\in\mathcal D_\alpha$ be fixed.  By
Proposition~\ref{LA} there is a point $\mu<0$ in the resolvent set
of $^\lambda\!\Delta$. Every $\mu<0$ is in the simply connected set
$\mathbb C\setminus\sigma_{ess}({^\lambda\!\Delta})$, cf.
Fig.~\ref{fig5}.  This  implies that
${^\lambda\!\Delta}-\mu:\lefteqn{\stackrel{\circ}{\phantom{\,\,>}}}H^2(\mathcal
G)\to L^2(\mathcal G)$ is an analytic  Fredholm operator function of
$\mu\in\mathbb C\setminus\sigma_{ess}({^\lambda\!\Delta})$. It is
known e.g.~\cite[Proposition A.8.4]{KozlovMaz`ya}, the spectrum  of
an analytic Fredholm operator function consists of isolated
eigenvalues of finite algebraic multiplicity. Thus $\sigma
({^\lambda\!\Delta})=\sigma_{ess}({^\lambda\!\Delta})\cup \sigma_d
({^\lambda\!\Delta})$.

{\it 5.} Recall that for all real  $\lambda\in \mathcal D_\alpha$
and $F,G\in L^2(\mathcal G)$
 the equality~\eqref{h6} is valid. To this equality we implement a variant of the Aguilar-Balslev-Combes argument. Namely, let $F$ and $G$ be in the set $\mathcal A$ of partial analytic vectors. Then by Proposition~\ref{p1}
$F\circ\vartheta_\lambda$ and $G\circ\vartheta_\lambda$ are
$L^2(\mathcal G)$-valued analytic functions of $\lambda$ in the disc
$\mathcal D_\alpha$.  In order to see that the equality~\eqref{h6}
extends by analyticity to all $\lambda\in\mathcal D_\alpha$, we
consider the resolvent $({^\lambda\!\Delta}-\mu)^{-1}$. By
Proposition~\ref{LA} $\lambda\mapsto {^\lambda\!\Delta}$ is an
analytic family in the sense of Kato, hence the resolvent is an
analytic function of the variables $\mu$ and $\lambda$ on the set
$\{(\mu,\lambda)\in\mathbb C\times\mathcal D_\alpha: \mu
\notin\sigma({^\lambda\!\Delta})\}$, e.g.~\cite[Theorem XII.7]{Simon
Reed iv}. Let $\mu$ be a negative number outside of the sector of
the m-sectorial operator ${^\lambda\!\Delta}$, cf.
Proposition~\ref{LA}. Then $\mu\notin\sigma({^\lambda\!\Delta})$ for
all $\lambda\in\mathcal D_\alpha$, and the resolvent $\mathcal
D_\alpha\ni\lambda\mapsto({^\lambda\!\Delta}-\mu)^{-1}$ is an
analytic function. As a consequence, the equality~\eqref{h6} extends
by analyticity to all $\lambda\in\mathcal D_\alpha$.

We take some $\lambda\in\mathcal D_\alpha$ with $\Im\lambda\neq 0$,
and consider the Left Hand Side (LHS) and the Right Hand Side (RHS)
of the equality~\eqref{h6} as functions of $\mu$.   The LHS is
meromorphic in $\mu\in\mathbb C\setminus\sigma_{ess}(\Delta)$ with
poles at the points of $\sigma_{d}(\Delta)$. Note that the spectrum
$\sigma(\Delta)$ is a subset of the half-line ${\mathbb R}_+$ as
$\Delta$ is a positive selfadjoint operator. On the other hand,  the
RHS is a meromorphic function  on the set $\mu\in\mathbb
C\setminus\sigma_{ess}({^\lambda\!\Delta})$. Therefore the RHS
provides the LHS with a meromorphic continuation in $\mu$ from
$\mathbb C\setminus\sigma_{ess}(\Delta)$ across the intervals
$(\nu_j,\nu_{j+1})$, $j\in\mathbb N$, between the thresholds to the
strips between the rays of $\sigma_{ess}({^\lambda\!\Delta})$, cf.
Fig.~\ref{fig5}.

 It is clear that the meromorphic continuation  can have poles only  at points  of $\sigma_d({^\lambda\!\Delta})$. Conversely, let $\mu_0\in\sigma_d({^\lambda\!\Delta})$, and let $\mathsf
P$ be the corresponding Riesz projection (i.e. the first order
residue of $({^\lambda\!\Delta}-\mu)^{-1}$ at the pole $\mu_0$). The
kernel $\ker({^\lambda\!\Delta}-\mu_0)\neq\{0\}$ is in the range of
$\mathsf P$.  Recall that the form $(\cdot,\cdot)_\lambda$ in
$L^2(\mathcal G)$ is nondegenerate, and by Proposition~\ref{p1} the
sets $\vartheta_\lambda[\mathcal A]$ and
$\vartheta_{\overline{\lambda}}[\mathcal A]$ are dense in
$L^2(\mathcal G)$. Therefore for some $F, G\in{\mathcal A}$ we must
have $ \bigl(\mathsf P
F\circ\vartheta_\lambda,G\circ\vartheta_{\overline{\lambda}}\bigr)_\lambda\neq
0$. Thus $\mu_0$ is a pole.

{\it 1.} The LHS of the equality~\eqref{h6} is an independent of $v$
analytic function of $\mu\in\Bbb C\setminus\sigma(\Delta)$. Hence
the meromorphic continuation of the LHS and its poles are
independent of the scaling function $v$. This together with the
assertion~{\it 5} implies that $\sigma_d({^\lambda\!\Delta})$ is
independent of  $v$. By the proven assertions~{\it 2} and~{\it 3}
the essential spectrum $\sigma_{ess}({^\lambda\!\Delta})$ is
independent of $v$ and $\sigma
({^\lambda\!\Delta})=\sigma_{ess}({^\lambda\!\Delta})\cup \sigma_d
({^\lambda\!\Delta})$.

{\it 4.} Let $\mu\in\sigma_d({^\lambda\!\Delta})$ be fixed. As
$\lambda$ changes continuously in the disk $\mathcal D_\alpha$ and
$\mu\notin\sigma_{ess}({^\lambda\!\Delta})$, the RHS of~\eqref{h6}
provides the LHS with one and the same meromorphic continuation to a
neighborhood of $\mu$. Therefore $\mu$ remains to be a pole of the
meromorphic continuation. Then  $\mu\in\sigma_d({^\lambda\!\Delta})$
by the assertion~{\it 5}.

{\it 6.}  Let $\lambda$ be a non-real number in the disk $\mathcal
D_\alpha$. Consider the projection
$$
\mathsf P=\operatorname{s-}\!\lim_{\epsilon\downarrow 0} i\epsilon
(\Delta-\mu_0-i\epsilon)^{-1}
$$
  onto the eigenspace of
the selfadjoint operator $\Delta$. Suppose that $\mu_0\in \mathbb
R\setminus\sigma({^\lambda \!\Delta})$ (then $\mu_0\neq \nu_j$ for
all $j\in\mathbb N$). Then for any $F,G\in\mathcal A$ the RHS
of~\eqref{h6} is an analytic function of $\mu$ in a neighborhood of
$\mu_0$. The equality~\eqref{h6} implies that $(\mathsf P F,G)=0$.
By Proposition~\ref{p1} the set $\mathcal A$ is dense in
$L^2(\mathcal G)$, and hence $\mathsf P=0$. Thus
$\ker(\Delta-\mu_0)=\{0\}$.

Now we assume that $\mu_0\in\mathbb R$ and
$\mu_0\in\sigma_d({^\lambda \!\Delta})$. Then the resolvent
$({^\lambda \!\Delta}-\mu)^{-1}$ has a pole at $\mu_0$. The sets
$\vartheta_\lambda[\mathcal A]$ and
$\vartheta_{\overline{\lambda}}[\mathcal A]$ are dense in
$L^2(\mathcal G)$. Hence there exist $F,G\in \mathcal A$, such that
$\mu_0$ is a pole for the RHS of~\eqref{h6}. The equality~\eqref{h6}
implies that $(\mathsf P F,G)\neq 0$, and thus
$\ker(\Delta-\mu_0)\neq\{0\}$.

{\it 7.} The RHS of~\eqref{h6} with $\Im\lambda> 0$, and therefore
the LHS, being defined on the dense subset $\mathcal A$ of
$L^2(\mathcal G)$,  has limits at the points $\Bbb
R\setminus\sigma({^\lambda\!\Delta})$
 as $\mu$ tends to the real line from $\mathbb C^+$. Since the set $\Bbb R\cap\sigma({^\lambda\!\Delta})$ is countable,  the Dirichlet Laplacian
$\Delta$ has no singular continuous spectrum,
e.g.~\cite[Theorem~XII.20]{Simon Reed iv}.$\Box$
\end{pf*}


\begin{thebibliography}{100}

\bibitem{AC} { J. Aguilar and J.M. Combes,} {\em A class of analytic perturbations for one-body
Schr¨odinger Hamiltonians,} Commun. Math. Phys. 22 (1971),
pp.~280-294



\bibitem{Aslan} { A. Aslanyan, L. Parnovski, and D. Vassiliev},
`Complex resonances in acoustic waveguides', {\em Q. J. Mech. Appl.
Math.} 53 (2000) 429-447.

\bibitem{BC}{ E. Balslev and J.M. Combes,} `Spectral properties of many-body Schr¨odinger
operators with dilatation analytic interactions', {\em Commun. Math.
Phys.} 22 (1971) 280–294


\bibitem{chr ST} { T. Christiansen,} `Scattering theory for manifolds with asymptotically cylindrical ends', {\em J. Funct. Anal.} 131 (1995) 499--530.

\bibitem{C}{ J. M. Combes,} `An algebraic approach to quantum scattering', unpublished.

\bibitem{Cycon} {H.~L. Cycon, R.~G. Froese, W. Kirsch, and B. Simon,}
{\em Schr\"{o}dinger operators, with application to quantum
mechanics and global geometry} (Springer-Verlag, New York, 1986).



\bibitem{DES} { P. Duclos, P. Exner, P. \v{S}\v{t}ov\'{\i}\v{c}ek}, `Curvature-induced resonances in a two dimensional Dirichlet tube', {\em Ann. Inst. H. Poincar\`{e}: Phys.Th\`{e}or.} 62 (1995)~81--101.


\bibitem{DEM} {P. Duclos, P.Exner, B. Meller,} `Exponential bounds on curvature-induced resonances in a two-dimensional Dirichlet tube', {\em Helv. Phys. Acta} 71 (1998)~133--162.

\bibitem{Edward}{ J. Edward,} `Eigenfunction decay and eigenvalue
accumulation for the Laplacian on asymptotically perturbed
waveguides', {\em J. London Math. Soc.}  59 (1999) 620--636.


\bibitem{FroHislop} { R. Froese and P. Hislop,} `Spectral analysis of second-order elliptic
operators on noncompact manifolds', {\em Duke Math. J.} 58 (1989)
103--129.

\bibitem{Guillope} { L. Guillop\'{e},} `Th\'{e}orie spectrale de quelques vari\'{e}tes \`{a} bouts', {\em Ann. Sci. \'{E}cole Norm. Sup.} 22 (1989) 137--160.




\bibitem{Hislop Sigal}{ P.~D. Hislop and I.~M. Sigal,} {\em Introduction to
spectral Theory: with applications to Schr\"{o}dinger Operators}
(Applied Mathematical Sciences 113, Springer-Verlag, 1996).

\bibitem{hermander}{ L. H\"{o}rmander,} {\em The analysis of linear partial differential operators III: Pseudodifferential operators} (Springer, 1985).




\bibitem{Hunziker} { W. Hunziker}, `Distortion analyticity and molecular resonance curves', {\em Ann. Inst. H. Poincare Phys. Theor.} 45 (1986) 339--358.


\bibitem{KalvinSiNum} { V. Kalvin}, `Perfectly matched layers for diffraction
gratings in inhomogeneous media. Stability and error estimates',
 {\em to appear in SINUM}.

\bibitem{Kato} { T. Kato,} {\em Perturbation theory for linear operators}
(Springer, Berlin-Heidelberg-New York, 1966).

\bibitem{KozlovMaz`ya} { V.~A. Kozlov and  V.~G. Maz'ya,} {\em Differential equations with
operator coefficients (with applications to boundary value problems
for partial differential equations)} (Berlin, Springer-Verlag,
1999).

\bibitem{KozlovMazyaRossmann}  { V.~A. Kozlov,  V.~G. Maz'ya, and J. Rossmann}, {\em Elliptic
boundary value problems in domains with point singularities}
(Mathematical Surveys and Monographs, vol. 52, American Mathematical
Society, 1997).



\bibitem{Lions Magenes} { J.-L. Lions and E. Magenes,} {\em Non-homogeneous boundary value
problems and applications I}  (Springer-Verlag, New York-Heidelberg,
1972).




\bibitem{MV1}{ R.~Mazzeo and A.~Vasy,} `Analytic continuation of the resolvent of the Laplacian on symmetric spaces of noncompact type', {\em
J. Funct. Anal.} 228 (2005) 311-368.

\bibitem{MV2}{ R.~Mazzeo and A.~Vasy,} `Scattering theory on SL(3)/SO(3): connections with quantum 3-body scattering',
{\em Proc. Lond. Math. Soc.}  94 (2007) 545-593.

\bibitem{Melrose} { R.~B. Melrose,} {\em The Atiyah-Patodi-Singer index theorem}  (Research Notes in Mathematics, Wellesley, 1993).


\bibitem{MelroseScat} { R.~B.  Melrose,} {\em  Geometric scattering theory} (Cambridge University Press, Cambridge, 1995).

\bibitem{Mueller} { W. M\"{u}ller and G. Salomonsen,}  `Scattering theory for the Laplacian on manifolds with bounded curvature', J. Funct. Anal. 253 (2007) 158--206.

\bibitem{Peetre} { J. Peetre,} `Another approach to elliptic boundary problems', {\em Comm. Pure Appl. Math.} 14 (1961) 711-731.

\bibitem{Simon Reed iv} { M. Reed and B. Simon,} {\em Metods of modern mathematical physics.
IV: Analysis of operators} (Academic Press, New York, 1972).

\bibitem{S} { J. Sj\"{o}strand,}
`Resonances for bottles and trace formulae', {\em Math. Nachr.} 221
(2001) 95-149.

\bibitem{SZ} { J. Sj\"{o}strand and M.~Zworski,} `The complex scaling method for scattering by strictly convex obstacles',
{\em Ark. Mat.} 33 (1995) 135-172.

\bibitem{Titchmarsh} { E.~C. Titchmarsh,} {\em
Introduction to the theory of Fourier's integrals} (Oxford Univ.
Press, 1937).

\bibitem{WZ}{ J.~Wunsch and M.~Zworski,} `Distribution of resonances for asymptotically Euclidean manifolds', {\em
J. Differ. Geom.} 55 (2000) 43-82.

\end{thebibliography}
\end{document}